\documentclass[11pt]{amsart}   	

\usepackage[margin=1in]{geometry}
									
\usepackage{amsmath,amssymb, amsfonts,amsthm}
\usepackage{hyperref}
\usepackage{mathtools}
\usepackage{tikz-cd}
\usetikzlibrary{decorations.markings}
\usepackage{enumitem}
\usepackage{makecell}
\setcellgapes{4pt}

\newtheorem{theorem}{Theorem}[section]
\newtheorem{lemma}[theorem]{Lemma}
\newtheorem{corollary}[theorem]{Corollary}
\newtheorem{proposition}[theorem]{Proposition}
\newtheorem{conjecture}[theorem]{Conjecture}

\theoremstyle{definition}
\newtheorem{definition}[theorem]{Definition}
\newtheorem{example}[theorem]{Example}

\newtheorem{remark}[theorem]{Remark}

\newcommand{\lperp}[1]{\prescript{\perp}{}{#1}}

\newcommand{\add}{\mathsf{add}}
\newcommand{\Filt}{\mathsf{Filt}}
\newcommand{\Gen}{\mathsf{Gen}}

\newcommand{\Cogen}{\mathsf{Cogen}}

\newcommand{\mods}{\mathsf{mod}}

\DeclareMathOperator{\Hom}{\mathrm{Hom}}
\DeclareMathOperator{\Ext}{\mathrm{Ext}}
\DeclareMathOperator{\im}{\mathrm{image}}
\DeclareMathOperator{\coker}{\mathrm{coker}}
\DeclareMathOperator{\End}{\mathrm{End}}

\newcommand{\T}{\mathcal T}
\newcommand{\F}{\mathcal F}
\newcommand{\W}{\mathcal W}
\newcommand{\V}{\mathcal V}
\newcommand{\E}{\mathcal E}
\newcommand{\J}{\mathcal J}
\newcommand{\Jinv}{
	\rotatebox[origin=c]{180}{$\J$}}
\renewcommand{\P}{\mathcal P}
\newcommand{\C}{\mathcal C}
\newcommand{\Db}{\mathcal{D}^b}
\newcommand{\rk}{\mathsf{rk}}

\usepackage{soul}

\title{$\tau$-perpendicular wide subcategories}

\author{Aslak Bakke Buan}
\author{Eric J. Hanson}
\thanks{This work was supported by grant number FRINAT 301375 from the Norwegian Research Council. The authors wish to thank Erlend D. B{\o}rve and H{\r a}vard U. Terland for many insightful conversations. They also wish to thank Haruhisa Enomoto for pointing out a mistake in the first version of this manuscript, and for sharing with them the examples from \cite{ringel} discussed in Remark~\ref{rem:notLattice?}.\\
\indent This is the accepted manuscript of a paper published in {\it Nagoya Mathematical Journal}. This version is available under a Creative Commons CC BY-NC-ND license. The official version of record is available at \url{https://doi.org/10.1017/nmj.2023.16}.  \copyright \ The Author(s)}
\address{Department of Mathematical Sciences, Norwegian University of Science and Technology (NTNU), 7491 Trondheim, NORWAY}
\email{aslak.buan@ntnu.no}
\email{ejhanso3@ncsu.edu}

\subjclass[2020]{
16G10, 16G20}

\keywords{$\tau$-tilting theory, perpendicular category, $\tau$-tilting reduction, $\tau$-cluster morphism categories}

\begin{document}
\maketitle

%%%%%%%%%%%%%%%%%%%%%%%%%%%%%%%%%%%%%%%%

\begin{abstract}
	Let $\Lambda$ be a finite-dimensional algebra. A wide subcategory of $\mathsf{mod}\Lambda$ is called \emph{left finite} if the smallest torsion class containing it is functorially finite. In this paper, we prove that the wide subcategories of $\mathsf{mod}\Lambda$ arising from $\tau$-tilting reduction are precisely the Serre subcategories of left finite wide subcategories. As a consequence, we show that the class of such subcategories is closed under further $\tau$-tilting reduction. This leads to a natural way to extend the definition of the ``$\tau$-cluster morphism category" of $\Lambda$ to arbitrary finite-dimensional algebras. This category was recently constructed by Buan--Marsh in the $\tau$-tilting finite case and by Igusa--Todorov in the hereditary case.
\end{abstract}

%%%%%%%%%%%%%%%%%%%%%%%%%%%%%%%%%%%%%%%%

\tableofcontents

%%%%%%%%%%%%%%%%%%%%%%%%%%%%%%%%%%%%%%%%

\section{Introduction}

In the study of module categories of rings and algebras, certain classes of 
subcategories play a prominent role.  {\em Torsion pairs} are pairs $(\T, \F)$ of subcategories, where the torsion classes $\T$
are characterized by being 
closed under extensions and factors, and $\F = \T^{\perp} \colon=\{X \mid \Hom(\T,X) = 0\}$.
Together with their triangulated siblings, the $t$-structures, such pairs are 
closely connected to classical tilting theory, e.g. via
the Brenner-Butler theorem \cite{BB}  and HRS-tilting \cite{HRS}.  
More recently, inspired by links to cluster combinatorics,  Adachi, Iyama and
Reiten defined support \linebreak $\tau$-tilting modules in \cite{AIR}.
They showed that
functorially finite torsion classes are
exactly those of the form $\Gen M$ (i.e. all modules which are generated by sums of
copies of $M$), where $M$ is a support $\tau$-tilting module. This
strengthens a classical result of Auslander and Smal{\o} \cite{AS}.
 
{\em Wide subcategories} are exact abelian subcategories. They were first considered
by Hovey \cite{hovey} in the setting of commutative noetherian rings.
The importance 
of such subcategories in dealing with categories $\mods \Lambda$ of finitely generated modules over finite-dimensional algebras has been highlighted by work
of Ingalls--Thomas \cite{IT} and Marks--{\v S}{\v t}ov\'i{\v c}ek \cite{MS}. 
In particular \cite{MS} shows that there is a natural injective map from the set of wide subcategories to the set of
torsion classes of $\mods \Lambda$, and also that there is an injective map from functorially finite torsion classes 
to functorially finite wide subcategories. The wide subcategories in the image of this map are called {\em left finite} wide subcategories, and
there is also a dual notion of {\em right finite} wide subcategories.
For $\tau$-tilting finite algebras, all wide subcategories and torsion classes are functorially finite, and
the above maps are actually bijections between finite sets.

Functorially finite wide subcategories are known to be exactly those which are equivalent to module categories,
and are hence of special interest. Examples of such are the left finite and right finite wide subcategories, and also the Serre subcategories. 
The latter are subcategories of $\mods \Lambda$ which are equivalent to $\mods (\Lambda/I)$ for $I$ generated
by an idempotent in $\Lambda$.
Another important source of functorially finite wide categories is the {\em $\tau$-perpendicular categories}, first considered by Jasso \cite{jasso}. These are categories given by $M^{\perp} \cap {^{\perp}(\tau M)} \cap P^{\perp}$, where $(M,P)$
is a pair of modules with $\Hom(M, \tau M) = 0$ and $P$ a projective module satisfying $\Hom(P,M) = 0$. 
These generalize both Serre subcategories and moreover classical Geigle-Lenzing perpendicular categories \cite{GL}, which have been much studied in the hereditary setting.
There is also a dual concept of $\tau^{-1}$-perpendicular categories. See Definition~\ref{def:tauInvPerp}.

Our first main result gives a characterization of $\tau$-perpendicular categories, showing how the different classes mentioned above are related.

\begin{theorem}[Theorem \ref{thm:mainA}]\label{thm:intro:mainA}
	Let $\Lambda$ be a finite-dimensional algebra and let $\W\subseteq \mods\Lambda$ be a wide subcategory. Then the following are equivalent.
	\begin{enumerate}
		\item $\W$ is a $\tau$-perpendicular subcategory of $\mods\Lambda$.
		\item $\W$ is a $\tau^{-1}$-perpendicular subcategory of $\mods\Lambda$.
		\item $\W$ is a Serre subcategory of a left finite wide subcategory.
		\item $\W$ is a Serre subcategory of a right finite wide subcategory.
		\item There exists a functorially finite torsion class $\T\subseteq\mods\Lambda$ and a functorially finite torsion-free class $\F\subseteq\mods\Lambda$ with $\T^\perp \subseteq \F$ such that $\W = \T\cap\F$.
	\end{enumerate}
\end{theorem}

We note that Serre subcategories of wide subcategories also occur in Asai and Pfeifer's classification of so-called ``wide intervals'' of torsion classes \cite{AP}. We discuss the relationship between Theorem~\ref{thm:intro:mainA}, the results of Asai and Pfeifer, and the ``brick labeling'' of the lattice of torsion classes in Remark~\ref{rem:AP}.

It is a consequence of Theorem~\ref{thm:intro:mainA} that the left-finite wide subcategories, right-finite wide subcategories, and Serre subcategories are examples of $\tau$-perpendicular subcategories. In particular, this leads to the following consequence.

\begin{corollary}[Corollary~\ref{cor:mainB}]\label{cor:intro:main}
		Let $\Lambda$ be a finite-dimensional algebra. Let $\V\subseteq \W \subseteq\mods\Lambda$ be a chain of subcategories such that $\V$ is a $\tau$-perpendicular subcategory of $\W$ and $\W$ is a $\tau$-perpendicular subcategory of $\mods\Lambda$. Then $\V$ is a $\tau$-perpendicular subcategory of $\mods\Lambda$.
\end{corollary}

Considering the finite poset $\mathcal{S}$ of all wide subcategories of a $\tau$-tilting finite algebra $\Lambda$, 
it was shown in \cite{BM_wide} that there is a natural definition of a category  $\mathfrak{W}(\Lambda)$,
with the elements in $\mathcal{S}$ as objects and maps parameterized by support $\tau$-rigid objects. Following \cite{HI_picture}, we call $\mathfrak{W}(\Lambda)$ the \emph{$\tau$-cluster morphism category} of $\Lambda$.
 The concept of (signed) $\tau$-exceptional sequences \cite{BM_exceptional} is closely related, as such sequences can be interpreted
 as compositions of irreducible maps in $\mathfrak{W}(\Lambda)$.
 This extended earlier work of Igusa-Todorov \cite{IT_exceptional}, who dealt with the hereditary case.
 The study of $\mathfrak{W}(\Lambda)$ was motivated by the link to the study of picture groups \cite{ITW} in the hereditary case, which was extended to the general 
 ($\tau$-tilting finite) case 
 in \cite{HI_picture}.
 
 As an application of Theorem \ref{thm:intro:mainA}, we show that one obtains a natural generalization of the above for all finite-dimensional algebras by restricting to  $\tau$-perpendicular subcategories. More precisely, we define a  category $\mathfrak{W}(\Lambda)$
 whose objects are the $\tau$-perpendicular subcategories of $\mods\Lambda$ and whose morphisms with source $\W$ are parameterized by the support $\tau$-rigid objects of $\W$. The following is then our second main result.
  
 \begin{theorem}[Theorem~\ref{thm:mainB}]\label{thm:intro:mainB}
 	Let $\Lambda$ be a finite-dimensional algebra. Then the $\tau$-cluster morphism category  $\mathfrak{W}(\Lambda)$ is a well-defined category.
 \end{theorem}
 
 The proof we give of Theorem \ref{thm:intro:mainB} also gives a significant simplification of the proof in the $\tau$-tilting finite case, given in \cite{BM_wide}.
 
 The paper is organized as follows. We first recall results and definitions concerning torsion pairs and
 $\tau$-tilting theory in Section \ref{sec:torsion}. Then we consider various classes of functorially finite 
 wide subcategories in Section \ref{sec:wide} and proceed by proving the first main theorem in Section \ref{sec:reductionCats}.
  We review a reduction formula for support $\tau$-rigid objects in Section \ref{sec:reductionFormulas}, which
 is used to prove the second main theorem in Section \ref{sec:catOfWideSubcats}. We conclude
 by working out a concrete example of a $\tau$-cluster morphism category in the final section.

\section{Torsion pairs and $\tau$-tilting theory}\label{sec:torsion}

In this section, we recall necessary background on torsion pairs and $\tau$-tilting theory.
Throughout this paper, $\Lambda$ shall always denote a finite-dimensional basic algebra over a field $K$, and $\mods \Lambda$ denotes
the category of finitely generated left $\Lambda$-modules. Furthermore, the Auslander-Reiten (AR) translate in
$\mods \Lambda$ is denoted by $\tau$.

The study of $\tau$-tilting theory has become instrumental in the study of finitely generated $\Lambda$-modules since its inception in \cite{AIR}. 
We follow the notation of \cite{BM_wide}, and denote $\C(\mods\Lambda):= \mods\Lambda \sqcup \mods\Lambda[1] \subseteq \Db(\mods\Lambda)$,
where $\Db(\mods\Lambda)$ denotes the bounded derived category of $\Lambda$.
An (usually assumed basic) object $U = M\sqcup P[1]\in \C(\mods\Lambda)$ is called a \emph{support $\tau$-rigid pair} if
\begin{enumerate}
	\item $M \in \mods\Lambda$ satisfies $\Hom(M,\tau M) = 0$.
	\item $P \in \mods\Lambda$ is projective and satisfies $\Hom(P,M) = 0$.
\end{enumerate}
If $U$ is basic, we denote by $\rk(U)$ the number of indecomposable direct summands of $U$ (up to isomorphism). If $\rk(U) = \rk(\Lambda)$, then $U$ is called  \emph{support $\tau$-tilting}. When $P = 0$, the module $U = M$ can be referred to as a $\tau$-rigid (or $\tau$-tilting if it is sincere and $\rk(U) = \rk(\Lambda)$) module.

By a subcategory of $\mods\Lambda$, we shall always mean a full subcategory which is closed under isomorphisms. Given such a subcategory $\mathcal{A}\subseteq \mods\Lambda$, we denote by $\P(\mathcal{A})$ the category of modules which are ext-projective in $\mathcal{A}$. That is, $Q \in \P(\mathcal{A})$ if and only if $\Ext^1(Q,X) = 0$ for all $X \in \mathcal{A}$.

Moreover, given a subcategory $\mathcal{A}\subseteq \mods\Lambda$, we denote by $\Gen(\mathcal{A})$ (resp. $\Cogen(\mathcal{A})$) the subcategory of $\mods\Lambda$ consisting of objects which are factors (resp. subobjects) of direct sums of objects in $\mathcal{A}$. We likewise denote by $\Filt(\mathcal{A})$ the subcategory of modules which admit finite filtrations whose subsequent subfactors all lie in $\mathcal{A}$. Given a module $X \in \mods\Lambda$, we define $\Gen X:= \Gen(\add X)$, etc., where $\add X$ is the subcategory of direct summands of finite direct sums of $X$. 

For any subcategory $\mathcal{A}\subseteq \mods\Lambda$, we associate two additional subcategories:
\begin{eqnarray*}
	\mathcal{A}^\perp &=& \{Y \in \mods\Lambda \mid \forall X \in \mathcal{A}: \Hom(X,Y) = 0\}\\
	\lperp{\mathcal{A}} &=& \{X \in \mods\Lambda \mid \forall Y \in \mathcal{A}: \Hom(X,Y) = 0\}.
\end{eqnarray*}
Given a module $X \in \mods\Lambda$, we likewise have $X^\perp := (\add X)^\perp$ and $\lperp{X} = \lperp{(\add X)}$.

Finally, we recall that a subcategory $\mathcal{A} \subseteq \mods\Lambda$ is called \emph{functorially finite} if for all $X \in \mods\Lambda$:
\begin{enumerate}
	\item There exists $A_X \in \mathcal{A}$ and $a_X:A_X\rightarrow X$ such that every morphism with source in $\mathcal{A}$ and target $X$ factors through $a_X$. The morphism $a_X$ is called a \emph{right $\mathcal{A}$-approximation}.
	\item There exists $A^X \in \mathcal{A}$ and $a^X:X\rightarrow A^X$ such that every morphism with source $X$ and target in $\mathcal{A}$ and target $X$ factors through $a^X$. The morphism $a^X$ is called a \emph{left $\mathcal{A}$-approximation}.
\end{enumerate}

We are now ready to discuss torsion pairs. A \emph{torsion pair} is a pair $(\T,\F)$ of subcategories of $\mods\Lambda$ such that $\T^\perp = \F$ and $\lperp{\F} = \T$. In this case, we call $\T$ a torsion class and $\F$ a torsion-free class. It is well known that a pair $(\T,\T^\perp)$ (resp. $(\lperp{\F},\F)$) is a torsion pair if and only if $\T$ is closed under extensions and quotients (resp. $\F$ is closed under extensions and subobjects). Moreover, given a torsion pair $(\T,\F)$, we have that $\T$ is functorially finite if and only if $\F$ is functorially finite \cite{smalo}.

If $(\T,\F)$ is a torsion pair, then every $M \in \mods\Lambda$ admits a unique exact sequence of the form
\begin{equation}\label{eqn:canonical}0\rightarrow t_\T(M)\xrightarrow{\iota} M\xrightarrow{q} f_\F(M)\rightarrow 0\end{equation}
with $t_\T(M) \in \T$ and $f_\F(M) \in \F$. In particular, the map $\iota$ is a minimal right $\T$-approximation and the map $q$ is a minimal left $\F$-approximation. We note that the operations $t_\T(-)$ and $f_{\F}(-)$ are both functorial.

We will need the following observation for
our discussion of Example~\ref{ex:notTauPerp}.

\begin{lemma}\label{lem:ffTorsionTorsionFree}
	Let $\T$ be a functorially finite torsion class and let $\F$ be a functorially finite torsion-free class. Then $\T\cap \F$ is functorially finite.
\end{lemma}
\begin{proof}
	We will show only that left $(\T \cap \F)$-approximations exist, as the argument for right approximations is analogous. Let $X \in \mods\Lambda$. Let $t^X: X\rightarrow T^X$ be a left $\T$-approximation of $X$ and let $f^X: T^X \rightarrow f_\F(T^X)$ be the left $\F$-approximation of $T^X$ coming from Equation~\ref{eqn:canonical}. We note that $f^X$ is surjective, and so $f_\F(T^X) \in \T\cap \F$. It is then straightforward to show that $f^X \circ t^X$ is a left $(\T\cap \F)$-approximation of~$X$.
\end{proof}

We now turn our attention to the well-established relationship between torsion pairs and support $\tau$-rigid objects.

Torsion pairs are closely related to support $\tau$-rigid objects, as the following shows.

\begin{theorem}\cite[Sections 2.2-2.3]{AIR}\label{thm:ftorsRigid}
	Let $\Lambda$ be a finite-dimensional algebra. Then
	\begin{enumerate}
		\item If $U = M\sqcup P[1] \in \C(\mods\Lambda)$ is support $\tau$-rigid, then both $\Gen M$ and $\lperp{(\tau M)} \cap P^\perp$ are functorially finite torsion classes in $\mods\Lambda$.
		\item If $U = M\sqcup P[1] \in \C(\mods\Lambda)$ is support $\tau$-tilting, then $\Gen M = \lperp{(\tau M)} \cap P^\perp$. Moreover, this association gives a bijection between support $\tau$-tilting objects in $\C(\mods\Lambda)$ and functorially finite torsion classes of $\mods\Lambda$.
		\item Let $\T \subseteq \mods\Lambda$ be a functorially finite torsion class, let $M \in \mods\Lambda$ be basic such that $\add M = \P(\T)$, and let $P \in \P(\mods\Lambda)$ be the maximal basic projective module which satisfies $\Hom(P,M) = 0$. Then $M\sqcup P[1]$ is support $\tau$-tilting and satisfies $\Gen M = \T = \lperp{(\tau M)}\cap P^\perp$.
	\end{enumerate}
\end{theorem}

Before continuing, we recall the following characterization of Auslander and Smal\o, which will be useful in several of our proofs.

\begin{proposition}\cite[Proposition~5.8]{AS} \label{prop:AStau}
	Let $M, N \in \mods\Lambda$. Then $\Hom(N,\tau M) = 0$ if and only if $\Ext^1(M,\Gen N) = 0$.
\end{proposition}

It is implicit in Theorem~\ref{thm:ftorsRigid} that any basic support $\tau$-rigid object is the direct summand of at least one support $\tau$-tilting object. In particular, we have the following.

\begin{theorem}\cite[Section~2.2]{AIR}\label{thm:bongartzDef}
	Let $U = M\sqcup P[1] \in \C(\mods\Lambda)$ be a basic support $\tau$-rigid object. Then
	\begin{enumerate}
		\item There exists a unique module $B_U \in \mods\Lambda$ such that $B_U\sqcup U$ is support $\tau$-tilting and $\add(B_U\sqcup M) = \P(\lperp{(\tau M)}\cap P^\perp)$. In particular, this means $$\Gen(B_U\sqcup M) = \lperp{(\tau M)\cap P^\perp} = \lperp{(\tau(B_U\sqcup M))}\cap P^\perp.$$
		\item There exists a unique object $C_U = N\sqcup Q[1] \in \C(\mods\Lambda)$ such that $C_U\sqcup U$ is support $\tau$-tilting and $\add(N\sqcup M) = \P(\Gen M)$. In particular, this means $$\Gen(N\sqcup M) = \Gen M =\lperp{(\tau(N\sqcup M))}\cap Q^\perp.$$
	\end{enumerate}
\end{theorem}

The module $B_U$ in Theorem~\ref{thm:bongartzDef} is called the \emph{Bongartz complement} of $U$. Following e.g. \cite{DIRRT,BM_exceptional}, we refer to $C_U$ in Theorem~\ref{thm:bongartzDef} as the \emph{co-Bongartz complement} of $U$.

\begin{remark}\
	\begin{enumerate}
		\item In \cite{AIR}, the Bongartz complement is only explicitly defined when $P = 0$ (so $U = M$ is a $\tau$-rigid module). Nevertheless, the more general definition is often given the same attribution. See e.g. \cite[Section~4]{DIRRT}. 
		\item If $P = 0$, then $B_U\sqcup M$ is a sincere module. In general, we can instead see $B_U\sqcup M$ as a sincere object in the Serre subcategory $P^\perp$. In this case, it is straightforward to show that $M$ is $\tau$-rigid in $P^\perp$ and that the Bongartz complement of $M$ in $P^\perp$ is precisely $B_U$. See \cite[Lemma~3.8]{BM_wide}.
	\end{enumerate}
\end{remark}

We next recall two results which give us a ``canonical decomposition'' of a support $\tau$-tilting pair. The first can be seen as a combination of \cite[Lemma~2.8]{IT} and \cite[Lemma~3.7]{MS}. We recall that a module $X$ in a subcategory $\mathcal{A}\subseteq \mods\Lambda$ is called split projective (in $\mathcal{A}$) if every epimorphism in $\mathcal{A}$ with target $X$ is split.

\begin{lemma}\label{lem:splitNonsplit}
	Let $\T \subseteq \mods\Lambda$ be a functorially finite torsion class and let $M\sqcup P[1] \in \C(\mods\Lambda)$ be the support $\tau$-tilting pair which satisfies $\P(\T) = \add M$. Then:
	\begin{enumerate}
		\item There is a decomposition $M = M_s \sqcup M_{ns}$ such that $M_s$ is split projective in $\T$ and no direct summand of $M_{ns}$ is split projective in $\T$. In particular, $\T = \Gen M_s$.
		\item Let
		$$\Lambda \xrightarrow{g}T_0\rightarrow T_1\rightarrow 0$$
		be an exact sequence such that $g$ is a minimal left $\T$-approximation. Then $\add T_0 = \add M_s$ and $\add T_1 = \add M_{ns}$.
	\end{enumerate}
\end{lemma}

The second result relates the direct summands $M_s$ and $M_{ns}$ to Bongartz and co-Bongartz complements.

\begin{proposition}\label{prop:bongartzSplitNonsplit}
	Consider the setup in Lemma~\ref{lem:splitNonsplit}. Then $M_s$ is the Bongartz complement of $M_{ns}\sqcup P[1]$ and $M_{ns}\sqcup P[1]$ is the co-Bongartz complement of $M_s$. In particular, we have
	$$\Gen M_s = \Gen M = \lperp{(\tau M)}\cap P^\perp = \lperp{(\tau M_{ns})}\cap P^\perp.$$
\end{proposition}

\begin{proof}
	It is shown in \cite[Lemma~4.5]{DIRRT} that $M_s$ is the Bongartz complement of $M_{ns}\sqcup P[1]$. Moreover, it is clear from Lemma~\ref{lem:splitNonsplit} that $\Gen M_s = \Gen M$. The result thus follows from Theorem~\ref{thm:bongartzDef}.
\end{proof}

We conclude this section with a brief description of the dual theory of $\tau^{-1}$-tilting. In order to state this in our context, given an indecomposable stalk complex $M[m] \in \Db(\mods\Lambda)$ with $M \in \mods\Lambda$, we denote
$$\overline{\tau}(M[m]) := \begin{cases} \tau M[m], & M \notin \P(\mods\Lambda)\\ \nu M[m-1] & M \in \P(\mods\Lambda),\end{cases}$$
where $\nu$ denotes the Nakayama functor. We then say a (usually assumed basic) object $U = I[-1]\sqcup M \in \C(\mods\Lambda)[-1]$ is \emph{support $\tau^{-1}$-rigid} if:
\begin{enumerate}
	\item $M \in \mods\Lambda$ and $\Hom(\tau^{-1} M, M) = 0$.
	\item $I \in \mods\Lambda$ is injective and $\Hom(M,I) = 0$.
\end{enumerate}
We likewise say $U$ is $\tau^{-1}$-tilting if $U$ is basic and $\rk(U) = \rk(\Lambda)$. It is shown in \cite[Section~2.2]{AIR} that $U \in \C(\mods\Lambda)[-1]$ is support $\tau^{-1}$-rigid (resp. support $\tau^{-1}$-tilting) if and only if there exists some support $\tau$-rigid (resp. support $\tau$-tilting) $V \in \C(\mods\Lambda)$ such that $U = \overline{\tau}V$. 

%%%%%%%%%%%%%%%%%%%%%%%%%%%%%%%

\section{Wide subcategories}\label{sec:wide}

In this section, we recall the definition and basic properties of wide subcategories and discuss important classes of examples of functorially finite wide subcategories: $\tau$-perpendicular subcategories, left/right finite wide subcategories, and Serre subcategories.

Recall that wide subcategories are exactly embedded abelian subcategories, and that a subcategory $\W\subseteq \mods\Lambda$ is wide if and only if it is closed under extensions, kernels, and cokernels.
It is well known that a wide subcategory $\W\subseteq\mods\Lambda$ is functorially finite if and only if it is equivalent to $\mods\Lambda_\W$ for some basic finite-dimensional algebra $\Lambda_\W $
(This is made explicit in \cite[Proposition~4.12]{enomoto}.) Given such a wide subcategory, we denote by $\rk(\W)$ the number of simple objects in $\W$ (or equivalently simple modules in $\mods\Lambda_\W$) up to isomorphism. We note that if $P \in \W$ is basic and $\P(\W) = \add(P)$, then $\rk(\W) = \rk(P)$. In particular, $\rk(\Lambda) = \rk(\mods\Lambda)$.

As wide subcategories are abelian categories in their own right, they have their own wide subcategories, torsion classes, and torsion-free classes. We will be concerned with such subcategories in the sequel, and so the following well-known and easily proved fact is useful.

\begin{proposition}\label{prop:ffReduction}
	Let $\W\subseteq\mods\Lambda$ be a wide subcategory.
	\begin{enumerate}
		\item Suppose that $\W$ is functorially finite and let $\mathcal{A}\subseteq \W$ be a functorially finite subcategory of $\W$. Then $\mathcal{A}$ is a functorially finite subcategory of $\mods\Lambda$.
		\item Let $\mathcal{V}$ be a wide subcategory of $\W$. Then $\mathcal{V}$ is a wide subcategory of $\mods\Lambda$.
	\end{enumerate}
\end{proposition}

We are now ready to define our main categories of interest.

\begin{definition}\label{def:jassoWide}
	A full subcategory $\W \subseteq \mods\Lambda$ is called a \emph{$\tau$-perpendicular subcategory} if there exists a support $\tau$-rigid object $U = M\sqcup P[1] \in \C(\mods\Lambda)$ such that $\W = \J(U)$, where
	\begin{equation}\label{eqn:jasso} \J(U):= (M\sqcup P)^\perp \cap \lperp{(\tau M)}.\end{equation}
\end{definition}

Such categories were first considered by Jasso \cite{jasso}, who proved
that they are equivalent to module categories, and hence they are
functorially finite. Actually, Jasso explicitly dealt with the case $P=0$,
but his proofs and statements can be easily modified. This is mentioned explicitly in \cite{DIRRT}, where it is also shown that such categories are
in fact wide. Summarizing, we have:

\begin{theorem}\cite[Theorem~3.8]{jasso}\cite[Theorems~4.12, 4.16]{DIRRT}\label{thm:rank}
	Let $U = M\sqcup P[1]\in \C(\mods\Lambda)$ be support $\tau$-rigid. Then $\J(U)$ is a functorially finite wide subcategory of $\mods\Lambda$. Moreover, if $M\sqcup P[1]$ is basic, then $\rk(\J(U)) + \rk(M) + \rk(P) = \rk(\Lambda)$.
\end{theorem}

By identifying a $\tau$-perpendicular subcategory $\W = \J(U)$ with a module category, one can consider the $\tau$-tilting theory of $\W$. That is, we consider the category $\C(\W) := \W\sqcup \W[1] \subseteq \C(\mods\Lambda)$. We then say an object $N\sqcup Q[1] \in \C(\W)$ is \emph{support $\tau$-rigid in $\W$} if:
\begin{enumerate}
	\item $N \in \W$ satisfies $\Hom(N,\tau_\W N) = 0$, where $\tau_W$ denotes the Auslander-Reiten (AR) translate in $\W$.
	\item $Q \in \W$ is projective in $\W$ and satisfies $\Hom(Q,N) = 0$.
\end{enumerate}
We emphasize that in general, objects which are projective in $\W$ may not be projective in $\mods\Lambda$. Likewise, we may have that $\tau_\W N \not \cong \tau N$, so in general we can have modules which
are {\em not} $\tau$-rigid in $\mods \Lambda$, but still are $\tau$-rigid
in $\W$. However, it is a direct consequence of Proposition \ref{prop:AStau} and the fact that $\W$ is an
exactly embedded subcategory that $\tau$-rigid (or projective) objects in $\mods \Lambda$ remain $\tau$-rigid 
(or projective) in $\W$.

Note also that since $\W$ is an exact subcategory of $\mods \Lambda$,
the category $\C(\W)$ can be considered as sitting either inside $\C(\mods\Lambda) \subseteq \Db(\mods\Lambda)$ (as we have defined it) or inside of $\Db(\W)$. Indeed, for $X, Y \in \W$, we have a canonical isomorphism
    $$\Hom_{\Db(\mods\Lambda)}(X,Y[1]) = \Hom_{\Db(\W)}(X,Y[1]).$$

Continuing in this way, if $V = N\sqcup Q[1]$ is $\tau$-rigid in $\W$, we denote
\begin{equation}\label{eqn:jassoSubcat}
	\J_\W(V) := (N\sqcup Q)^\perp \cap \lperp{(\tau_W N)} \cap \W. \end{equation}

We also have the following dual concept. 

\begin{definition}\label{def:tauInvPerp}
	A subcategory $\W \subseteq \mods\Lambda$ is called a \emph{$\tau^{-1}$-perpendicular subcategory} if there exists a support $\tau^{-1}$-rigid object $U = I[-1]\sqcup M \in \C(\mods\Lambda)[-1]$ such that
	$$\W = \Jinv(U):= \lperp{(M\sqcup I)} \cap (\tau^{-1} M)^\perp.$$
\end{definition}

We show as part of Theorem~\ref{thm:intro:mainA} that $\tau^{-1}$-perpendicular subcategories and  $\tau$-perpendicular subcategories coincide. Moreover, it will be a consequence of Theorem~\ref{thm:intro:mainA} that not all functorially finite wide subcategories are $\tau$-perpendicular subcategories. See Example~\ref{ex:notTauPerp}.

We proceed to discuss another central class of functorially finite wide subcategories.
These arise from applying the so-called Ingalls--Thomas bijections \cite{IT,MS} to functorially finite torsion classes and torsion-free classes. We recall these constructions now.

\begin{definition}\label{def:IT}\
	\begin{enumerate}
		\item Let $\T\subseteq\mods\Lambda$ be a torsion class. The \emph{left wide subcategory} of $\mods\Lambda$ corresponding to $\T$ is
		$$\W_L(\T) := \{X \in \T \mid (Y\in\T,f:Y\rightarrow X)\implies \ker f \in \T\}.$$
		\item Let $\F\subseteq \mods\Lambda$ be a torsion-free class. The \emph{right wide subcategory} of $\mods\Lambda$ corresponding to $\F$ is 
		$$\W_R(\F) := \{X \in \F \mid (Y\in\F,f:X\rightarrow Y)\implies \coker f \in \F\}.$$
	\end{enumerate}
\end{definition}

One of the key results of \cite{IT} (hereditary case) and \cite{MS} (general case) is that for any wide subcategory $\W$, one has
$$\W_L(\Filt\Gen(\W)) = \W = \W_R(\Filt\Cogen(\W)).$$
They also show that $\Filt\Gen(W_L(\T)) = \T$ (resp. $\Filt\Cogen\W_R(\F) =\F$), when $\T$ (resp. $\F$) is
functorially finite.
Following Asai \cite{asai_semibricks}, a wide subcategory $\W \subseteq\mods\Lambda$ is called \emph{left finite} (resp. \emph{right finite}) if it is of the form $\W_L(\T)$ (resp. $\W_R(\F)$) for some functorially finite torsion class $\T$ (resp. torsion-free class $\F$). It is straightforward that if $\W$ is either left finite or right finite, then it is functorially finite. The converse, however, does not hold in general. See \cite[Example~3.13]{asai_semibricks} or Example~\ref{ex:notTauPerp}.

We conclude this section by discussing a well understood class of functorially finite wide subcategories, namely the
Serre subcategories. A subcategory $\mathcal{S}$ is \emph{Serre} if for any short exact sequence
	$$0\rightarrow X\rightarrow Y \rightarrow Z\rightarrow 0$$
	in $\mods\Lambda$, we have $Y \in \mathcal{S}$ if and only if $X,Z \in \mathcal{S}$. That is, $\mathcal{S}$ is closed under extensions, quotients, and subobjects.

Serre subcategories are indeed examples of wide subcategories. In fact, they are also both torsion classes and torsion-free classes, as the following shows. 

\begin{proposition}\label{prop:serre}
	Let $\mathcal{S} \subseteq \mods\Lambda$ be a subcategory. Then the following are equivalent:
	\begin{enumerate}
		\item $\mathcal{S}$ is a Serre subcategory.
		\item $\mathcal{S}$ is any two of a torsion class, a torsion-free class, and a wide subcategory.
		\item $\mathcal{S}$ is a torsion class, a torsion-free class, and a wide subcategory.
		\item $\mathcal{S} = P^\perp$ for some projective $P \in \P(\mods\Lambda)$.
		\item $\mathcal{S}$ is a wide subcategory and every object which is simple in $\mathcal{S}$ is simple in $\mods\Lambda$.
	\end{enumerate}
\end{proposition}

\begin{proof}
    The equivalence of (1), (4), and (5) is contained in \cite[Proposition~5.3]{GL}, and the equivalence of (1), (2), and (3) follows straightforwardly from the definitions.
\end{proof}

As useful consequences, we obtain the following corollaries.

\begin{corollary}\label{cor:serreProj}
	There is a bijection between isomorphism classes of basic projective modules in $\P(\mods\Lambda)$
and Serre subcategories of $\mods\Lambda$ given by $P\mapsto P^\perp$.
\end{corollary}

\begin{corollary}\label{cor:serreFinite}
	Let $\mathcal{S}\subseteq \mods\Lambda$ be a Serre subcategory. Then $\mathcal{S}$ is both a left finite wide subcategory and a right finite wide subcategory.
\end{corollary}

\begin{proof}
	Note that $\mathcal{S}$ is a wide subcategory which satisfies $\W_L(\mathcal{S}) = \mathcal{S} = \W_R(\mathcal{S})$ by Proposition~\ref{prop:serre}. The result then follows from Theorem \ref{thm:ftorsRigid}(1) and item (4) of Proposition~\ref{prop:serre}.
\end{proof}

%%%%%%%%%%%%%%%%%%%%%%%%%%%%%%%%%%
\section{Characterizing $\tau$-perpendicular subcategories}
\label{sec:reductionCats}

In this section, we give the proof of Theorem~\ref{thm:intro:mainA}, restated as Theorem~\ref{thm:mainA} below. This characterizes $\tau$-perpendicular subcategories as precisely the Serre subcategories of left-finite and right-finite wide subcategories of $\mods\Lambda$.

The following technical result will be useful for the proof.

\begin{proposition}\cite[Proposition~5.2.1]{BTZ}\label{prop:minExtending}
	Let $(\T,\F)$ be a torsion pair in $\mods\Lambda$ and let $X \in \W_R(\F)$. Then $X$ is simple in $\W_R(\F)$ if and only if the following hold:
	\begin{enumerate}
		\item Every proper factor of $X$ lies in $\T$.
		\item If $0\rightarrow X \rightarrow Y\rightarrow Z\rightarrow 0$ is a nonsplit exact sequence and $Z \in \T$, then $Y \in \T$.
		\item $X \in \F$.
	\end{enumerate}
\end{proposition}

\begin{remark}\
\begin{enumerate}
	\item Our statement of Proposition~\ref{prop:minExtending} is actually the dual of \cite[Proposition~5.2.1]{BTZ}.
	\item The simple objects of $\W_R(\F)$ are given an alternative characterization in terms of 2-term simple-minded collections in \cite{asai_semibricks}. The characterization in \cite{BTZ}, on the other hand, shows that the simple objects of $\W_R(\F)$ are precisely the ``minimal extending modules" for the torsion class $\T$, introduced in \cite{BCZ}.
\end{enumerate}
\end{remark}

We now start building towards our proof of Theorem~\ref{thm:intro:mainA} with the following lemmas.

\begin{lemma}\label{lem:finiteIsJasso} Let $(\T,\F)$ be a functorially finite torsion pair and let $M\sqcup P[1]$ be the support $\tau$-tilting object  in $\C(\mods\Lambda)$ for which $\add(M) = \P(\T)$. Write $M \cong M_s\sqcup M_{ns}$ as in Lemma~\ref{lem:splitNonsplit}. Then:
	\begin{enumerate}
		\item $\W_L(\T) = \J(M_{ns}\sqcup P[1])$.
		\item $\W_R(\F) = \J(M_s)$.
	\end{enumerate}
In particular, any wide subcategory of $\mods\Lambda$ which is either left finite or right finite is also a $\tau$-perpendicular subcategory.
\end{lemma}

We note that (1) also appears as \cite[Equation~1.2]{yurikusa}.

\begin{proof}
	(1) It is shown in \cite[Lemma~3.8]{MS} that $\W_L(\T) = M_{ns}^\perp \cap \lperp{(\tau M)} \cap P^\perp$. Moreover, by Proposition~\ref{prop:bongartzSplitNonsplit}, we have that $\lperp{(\tau M)} \cap P^\perp = \lperp{(\tau M_{ns})} \cap P^\perp$. This proves the result.
	
	(2) First let $X$ be a simple object of $\W_R(\F)$. (Note that $X$ is not necessarily simple in $\mods\Lambda$.) We will show that $X \in \J(M_s)$. Since $\J(M_s)$ is closed under extensions, this will imply that $\W_R(\F) \subseteq \J(M_s)$.
	
	We first note that $\W_R(\F) \subseteq \F = M_s^\perp$, so we need only show that $\Hom(X,\tau M_s) = 0$. Suppose to the contrary that $\Hom(X,\tau M_s) \neq 0$. By Proposition~\ref{prop:AStau} this means there exists $X' \in \Gen X$ and a nonsplit exact sequence of the form
	$$0 \rightarrow X' \rightarrow E\rightarrow M_s\rightarrow 0.$$
	By Proposition~\ref{prop:minExtending}, we note that $X'$ cannot be a proper quotient of $X$. Indeed, if this were the case, we would have $X' \in \T \subseteq \lperp{(\tau M)}$, a contradiction. Therefore, we can assume that $X' = X$. Applying Proposition~\ref{prop:minExtending} again, this implies that $E \in \T$. Since $M_s$ is split projective in $\T$, this is a contradiction.
	
	 Now let $Y \in \J(M_s)$. It is clear that $Y \in \F = M_s^\perp$. Thus let $Z \in \F = M_s^\perp$ and $g:Y\rightarrow Z$. We then have an exact sequence
	$$0 = \Hom(M_s,Z)\rightarrow \Hom(M_s,\coker g)\rightarrow \Ext^1(M_s,\im g) = 0,$$
	where the last term is zero by Proposition~\ref{prop:AStau} and the fact that $\im g$ is a quotient of $Y$. We conclude that $\coker g \in \F$, and therefore $Y \in \W_R(\F)$. This completes the proof.
\end{proof}

\begin{lemma}\label{lem:serreReduction2}
Let $\V \subseteq \W \subseteq \mods\Lambda$ be a chain of subcategories such that $\V$ is a Serre subcategory of $\W$ and $\W$ is a $\tau$-perpendicular subcategory of $\mods\Lambda$. Then $\V$ is a $\tau$-perpendicular subcategory of $\mods\Lambda$.
\end{lemma}

\begin{proof}
    Let $U = M\sqcup P[1] \in \C(\mods\Lambda)$ be support $\tau$-rigid and let $\mathcal{S}$ be a Serre subcategory of $\J(U)$.
	 By Proposition~\ref{prop:serre}, there exists $Q \in \P(\J(U))$ so that $\mathcal{S} = Q^\perp \cap \J(U)$.
	It follows from \cite[Proposition 3.14]{jasso} and Theorem \ref{thm:ftorsRigid} that
	$Q = f_{M^{\perp}}(B)$ for some direct summand $B$ of the Bongartz complement $B_U$ of $U$.
We then have an exact sequence $$ 0 \to t_{\Gen M}(B) \to B \to Q \to 0,$$
and since  $t_{\Gen M}(B)$ is in $\Gen M$, we have $\Hom(t_{\Gen M}(B), \J(U)) = 0$.
Hence we have $\mathcal{S} = Q^{\perp} \cap \J(U) = B^\perp \cap \J(U) = \J(B\sqcup U)$, and so $\mathcal{S}$ is a $\tau$-perpendicular subcategory of $\mods\Lambda$.
\end{proof}

We are now ready to prove our first main result.

\begin{theorem}[Theorem~\ref{thm:intro:mainA}]\label{thm:mainA}
	Let $\Lambda$ be a finite-dimensional algebra and let $\W\subseteq \mods\Lambda$ be a wide subcategory. Then the following are equivalent.
	\begin{enumerate}
		\item $\W$ is a $\tau$-perpendicular subcategory of $\mods\Lambda$.
		\item $\W$ is a $\tau^{-1}$-perpendicular subcategory of $\mods\Lambda$.
		\item $\W$ is a Serre subcategory of a left finite wide subcategory.
		\item $\W$ is a Serre subcategory of a right finite wide subcategory.
		\item There exists a functorially finite torsion class $\T\subseteq\mods\Lambda$ and a functorially finite torsion-free class $\F\subseteq\mods\Lambda$ with $\T^\perp \subseteq \F$ such that $\W = \T\cap\F$.
	\end{enumerate}
\end{theorem}

\begin{proof}
	$(1\iff 2):$ Recall that if $U \in \C(\mods\Lambda)$ is support $\tau$-rigid, then $\overline{\tau}U \in \C(\mods\Lambda)[-1]$ is support $\tau^{-1}$-rigid, and moreover that every support $\tau^{-1}$-rigid object in $\C(\mods\Lambda)[-1]$ occurs in this way. Thus suppose $U = M\sqcup P[1]$ is support $\tau$-rigid and write $M \cong M_p\sqcup M_{np}$, where $M_p \in \P(\mods\Lambda)$ and $M_{np}$ has no projective direct summand. We note that $\lperp{(\nu M_p)} = M_p^\perp$ and $\lperp{(\nu P)} = P^\perp$. Moreover, we have $\overline{\tau}(U) = (\nu M_p)[-1]\sqcup \tau M_{np} \sqcup \nu P$, where $\tau M_{np}$ has no injective direct summand. This means
	\begin{eqnarray*}
		\J(U) &=& M_{np}^\perp \cap \lperp{(\tau M_{np})} \cap (M_p\sqcup P)^\perp\\
		&=& (\tau^{-1}\tau M_{np})^\perp \cap \lperp{(\tau M_{np})} \cap \lperp{(\nu M_p \sqcup \nu P)}\\
		&=& \Jinv(\overline{\tau}(U)).
	\end{eqnarray*}
	This proves the result.
	
	$(1\implies 5)$: This follows from Theorem~\ref{thm:ftorsRigid} and the definition of $\J(U)$.
	
	$(5\implies 3)$: Write $\W = \T\cap \F$ with $\T$ a functorially finite torsion class and $\F$ a functorially finite torsion-free class. We will first show that $\W\subseteq \W_L(\T)$ using an argument similar to \cite[Lemma~3.8]{MS}. Let $X \in \W$ and let $g:Y\rightarrow X$ be a morphism in $\T$. Note that $\im g \in \W$ since it is a subobject of $X$ and a quotient of $Y$. Now consider the canonical exact sequence with respect to the torsion pair $(\lperp{\F},\F)$:
	$$0\rightarrow t_{(\lperp{\F})}(Y) \rightarrow Y \rightarrow f_\F (Y)\rightarrow 0.$$
	By assumption, we have $f_\F(Y) \in \T\cap \F = \W$ and $t_{(\lperp{\F})}(Y) \in \lperp{\F} \subseteq \T$. In particular, we have that $\Hom(t_{(\lperp{\F})}(Y), \im g) = 0$, and so $\im g$ is a quotient of $f_\F(Y)$. That is, we obtain the following diagram with rows and columns exact:
	$$\begin{tikzcd}
		&&0\arrow[d]&0\arrow[d]\\
		0\arrow[r] &t_{(\lperp{\F})}(Y) \arrow[d,equals]\arrow[r] &\ker g\arrow[d]\arrow[r]&\ker g/t_{(\lperp{\F})}(Y)\arrow[d]\arrow[r]&0\\
		0\arrow[r] &t_{(\lperp{\F})}(Y)\arrow[r]&Y\arrow[d]\arrow[r]&f_\F(Y)\arrow[d]\arrow[r]&0\\
		&&\im g \arrow[d]\arrow[r,equals]&\im g \arrow[d]\\
		&&0&0
	\end{tikzcd}$$
Now $\ker g/t_{(\lperp{\F})} \in \W\subseteq \T$ since $\W$ is closed under kernels. Therefore $\ker g \in \T$ since $\T$ is closed under extensions. We conclude that $X \in \W_L(\T)$.

We will now show that $\W$ is a Serre subcategory of $\W_L(\T)$. Let
$$0\rightarrow X\rightarrow Y\rightarrow Z\rightarrow 0$$
be a short exact sequence in $\W_L(\T)$. It is clear that if $X,Z \in \W$ then $Y \in \W$. Thus suppose $Y \in \W = \T\cap \F$. It follows that $X \in \W$ since it is in $\W_L(\T) \subseteq \T$ and $\F$ is closed under subobjects. Since $\W$ is closed under cokernels, it follows that $Z \in \W$ as well.
	
	$(5\implies 4)$: The proof is dual to that of $(5\implies 3)$, but we include the details here for convenience. We will first show that $\W \subseteq \W_R(\F)$. Let $X \in \W$ and let $g:X\rightarrow Y$ be a morphism in $\F$. We note that $\im g \in \W$ since it is a quotient of $X$ and a subobject of $Y$. Now consider the canonical exact sequence with respect to the torsion pair $(\T,\T^\perp)$:
	$$0 \rightarrow t_\T(Y) \rightarrow Y \rightarrow f_{(\T^\perp)}(Y)\rightarrow 0.$$
	By assumption, we have $t_\T(Y) \in \T\cap \F = \W$ and $f_{(\T^\perp)}(Y) \in \T^\perp \subseteq \F$. In particular, we have $\im g \subseteq t_\T(T)$. Therefore, we have an exact sequence
	$$0\rightarrow t_\T(Y)/\im g \rightarrow \coker f \rightarrow f_{(\T^\perp)}(Y)\rightarrow 0.$$
	Since $\W$ is wide, we have $t_\T(Y)/\im g \in \W\subseteq \F$. Since $\F$ is closed under extensions, this implies that $\coker g \in \F$. We conclude that $X \in \W_R(\F)$.
	
We will now show that $\W$ is a Serre subcategory of $\W_R(\F)$. Let
$$0\rightarrow X\rightarrow Y\rightarrow Z\rightarrow 0$$
be a short exact sequence in $\W_R(\F)$. It is clear that if $X,Z \in \W$ then $Y \in \W$. Thus suppose $Y \in \W = \T\cap \F$. It follows that $Z \in \W$ since it is in $\W_R(\F) \subseteq \F$ and $\T$ is closed under subobjects. Since $\W$ is closed under kernels, it follows that $X \in \W$ as well.
	
	$(3\implies 1)$: Let $\W\subseteq\W_L(\T)$ be a Serre subcategory of a left finite wide subcategory of $\mods\Lambda$. By Lemma~\ref{lem:finiteIsJasso}, we have that $\W_L(\T)$ is $\tau$-perpendicular in $\mods\Lambda$. It then follows from Lemma~\ref{lem:serreReduction2} that $\W$ is $\tau$-perpendicular in $\mods\Lambda$ as well.
	
	$(4\implies 1)$: The proof is analogous to that of $(3\implies 1)$.
\end{proof}

\begin{remark}\label{rem:AP}
    We note that the equivalences between (3), (4), and (5) in Theorem~\ref{thm:mainA} can also be deduced from \cite[Corollary~6.8]{AP} in Asai and Pfeifer's work on ``wide intervals'' in the lattice of torsion classes. (They deduce this corollary after working with wide subcategories and torsion classes which are not necessarily functorially finite.) This yields a characterization of $\tau$-perpendicular subcategories using the ``brick labeling'' of the lattice of torsion classes, defined in \cite{asai_semibricks} for functorially finite torsion classes and in \cite{BCZ,DIRRT} for all torsion classes. Namely, let $\mathcal{T}$ be a functorially finite torsion class, and choose a set of cover relations of the form $\mathcal{U}_i \subsetneq \mathcal{T}$. Then $\left[ \bigcap_i \mathcal{U}_i,\mathcal{T}\right]$ is a ``wide interval'' by \cite[Theorem~5.2]{AP}. By definition, this means that the intersection $\mathcal{W} := \left(\bigcap_i \mathcal{U}_i\right)^\perp \cap \mathcal{T}$ is a wide subcategory. Now the torsion class $\bigcap_i \mathcal{U}_i$ is functorially finite by \cite[Theorem~3.14]{jasso} (see also \cite[Corollary~6.8]{AP}), and so the intersection $\mathcal{W}$ is a $\tau$-perpendicular subcategory by Theorem~\ref{thm:mainA}. Moreover, the simple objects of $\mathcal{W}$ are precisely the ``brick labels'' of the chosen cover relations by \cite[Theorem~4.2(3)]{AP}. Finally, all $\tau$-perpendicular subcategories will be of this form, again by \cite[Corollary~6.8]{AP} and Theorem~\ref{thm:mainA}.
\end{remark}

\begin{remark}
\sloppy Due to the equivalence between $\tau$-perpendicular and  $\tau^{-1}$-perpendicular subcategories of $\mods\Lambda$, we will dispense with discussing the support $\tau^{-1}$-rigid case for the remainder of this paper. We nevertheless remark that the majority of our results can be restated for \linebreak $\tau^{-1}$-perpendicular subcategories by applying Theorem~\ref{thm:mainA}.
\end{remark}

We conclude this section by tabulating several consequences of Theorem~\ref{thm:mainA} and the preceding lemmas.

\begin{corollary}\label{cor:corank1}
	Let $\W \subseteq \mods\Lambda$ be a functorially finite wide subcategory with $\rk(\W) + 1 = \rk(\Lambda)$. Then the following are equivalent.
	\begin{enumerate}
	    \item $\W$ is a left finite wide subcategory.
	    \item $\W$ is a right finite wide subcategory.
	    \item $\W$ is a $\tau$-perpendicular subcategory.
	\end{enumerate}
\end{corollary}

\begin{proof}
    The implications $(1\implies 3)$ and $(2\implies 3)$ are contained in Lemma~\ref{lem:finiteIsJasso}. Thus assume (3).
    By Theorem~\ref{thm:mainA}, there exists a left finite wide subcategory $\V \subseteq\mods\Lambda$ such that $\W$ is a Serre subcategory of $\V$. Lemma~\ref{lem:finiteIsJasso} and Theorem~\ref{thm:rank} then imply that
    $$\rk(\Lambda) - 1 = \rk(\W) \leq \rk(\V) \leq \rk(\Lambda).$$
    Now if $\rk(\W) = \rk(\V)$, then $\W = \V$ as a consequence of the same lemma and theorem. In particular, $\W$ is a left finite wide subcategory of $\mods\Lambda$ in this case. Otherwise, $\V = \mods\Lambda$ by the same argument, and so $\W$ is a left finite wide subcategory of $\mods\Lambda$ by Corollary~\ref{cor:serreFinite}. We conclude that (3) implies (1). The proof that (3) implies (2) is identical.
\end{proof}

In \cite[Example~3.13]{asai_semibricks}, Asai gives an example of a functorially finite wide subcategory which is right finite but not left finite. (It is, however, a Serre subcategory of a left finite wide subcategory, consistent with Theorem~\ref{thm:mainA}.) By modifying Asai's example, we obtain an example of a functorially finite wide subcategory which is not a $\tau$-perpendicular subcategory.

\begin{example}\label{ex:notTauPerp}
    Consider the quiver
    $$\begin{tikzcd}[column sep = 0.75cm]Q=&
1\arrow[r,yshift=0.1cm,"\alpha_1"above]\arrow[r,yshift=-0.1cm,"\alpha_2" below]&2\arrow[r,yshift=0.1cm,"\beta_1"above]\arrow[r,yshift=-0.1cm,"\beta_2" below]&3
    \end{tikzcd}$$
    and let $\Lambda = KQ/(\beta_2\alpha_1,\beta_1\alpha_2)$. Consider the $\Lambda$-modules
    $$\begin{tikzcd}[column sep = 0.75cm]X_1=&
K\arrow[r,yshift=0.1cm,"1"above]\arrow[r,yshift=-0.1cm,"0" below]&K\arrow[r,yshift=0.1cm,"1"above]\arrow[r,yshift=-0.1cm,"0" below]&K, &&X_2=&
K\arrow[r,yshift=0.1cm,"0"above]\arrow[r,yshift=-0.1cm,"1" below]&K\arrow[r,yshift=0.1cm,"0"above]\arrow[r,yshift=-0.1cm,"1" below]&K.
    \end{tikzcd}$$
    We will demonstrate that $\W:= \Filt(X_1,X_2)$ is a functorially finite wide subcategory of $\mods\Lambda$ which is not a $\tau$-perpendicular subcategory.
    
    We first note that $\Hom(X_1,X_2) = 0 = \Hom(X_2,X_1)$ and that $\End(X_1) \cong K \cong \End(X_2)$. Moreover, it is straightforward to show that
    $$\begin{tikzcd}[column sep = 0.75cm]\tau X_1=&
K\arrow[r,yshift=0.1cm,"1"above]\arrow[r,yshift=-0.1cm,"0" below]&K\arrow[r,yshift=0.1cm,"0"above]\arrow[r,yshift=-0.1cm,"0" below]&0, &&\tau X_2=&
K\arrow[r,yshift=0.1cm,"0"above]\arrow[r,yshift=-0.1cm,"1" below]&K\arrow[r,yshift=0.1cm,"0"above]\arrow[r,yshift=-0.1cm,"0" below]&0.
    \end{tikzcd}$$
    In particular, $\Hom(X_1,\tau X_2) = 0 = \Hom(X_2,\tau X_1)$ and every morphism $X_1\rightarrow \tau X_1$ (or $X_2\rightarrow \tau X_2$) factors through the injective $I(2)$.
    By the Auslander-Reiten formulas, we conclude that $\Ext^1(X_i,X_j)$ for any $i,j \in \{1,2\}$. This means $\W$ is a wide subcategory equivalent to the module category of a semisimple algebra, so in particular $\W$ is functorially finite.
    
    We will now show that $\W$ is not a $\tau$-perpendicular subcategory. By Corollary~\ref{cor:corank1}, the fact that $\rk(\W) = \rk(\Lambda)-1$ means we need only show that $\W$ is not left finite. To see this, we note that the Serre subcategory $P(3)^\perp$ is equivalent to the module category of the Kronecker path algebra. Thus we can consider the Ingalls-Thomas bijection (Definition~\ref{def:IT}) $\W_L^{P(3)^\perp}$ in the category $P(3)^\perp$. Then the wide subcategory $\V = \W_L^{P(3)^\perp}(\Filt\Gen(\W)\cap P(3)^\perp) = \W_L^{P(3)^\perp}(\Filt\Gen(\tau X_1,\tau X_2))$ is the additive closure of a pair of homogeneous tubes having $\tau X_1$ and $\tau X_2$ at their mouths. It is then a well-known fact that right $\V$-approximations will in general not exist, and
    hence $\V$ is not functorially finite (in either $P(3)^\perp$ or in $\mods\Lambda$). In particular, this means $\V$ is not a left finite wide subcategory of $P(3)^\perp$, and so $\Filt\Gen(\W)\cap P(3)^\perp$ is not functorially finite. Since $P(3)^\perp$ is a functorially finite torsion-free class of $\mods\Lambda$ (see Proposition~\ref{prop:serre} and Corollary~\ref{cor:serreFinite}), 
  Lemma~\ref{lem:ffTorsionTorsionFree} then implies that the torsion class $\Filt\Gen(\W)$ is not functorially finite; that is, that $\W$ is not left finite. (Note that, as a consequence of Corollary~\ref{cor:corank1}, $\W$ is not right finite either. This can also be seen directly using duality.)
\end{example}

\begin{remark}\label{rem:inclusions}
Summarizing, we have the following inclusions of classes of subcategories of $\mods\Lambda$:
	$$\begin{tikzcd}[row sep = 20pt]
	& \{\text{wide} \}& \\
	& \{\text{functorially finite wide} \} \arrow[draw=none]{u}[sloped,auto=false]{\text{\LARGE{$\subseteq$}}}\arrow[draw=none,u,xshift=0.5cm,"(1)"right]& \\
 & \{\text{$\tau$-perpendicular}\} \arrow[draw=none]{u}[sloped,auto=false]{\text{\LARGE{$\subseteq$}}}\arrow[draw=none,u,xshift=0.5cm,"(2)"right] & \\
 \{\text{left finite wide}\} 
 \arrow[draw=none]{ur}[sloped,auto=false]{\text{\LARGE{$\subseteq$}}}\arrow[draw=none,ur,xshift=-0.5cm,"(3)"left]& & \{\text{right finite wide}\} 
 \arrow[draw=none]{ul}[sloped,auto=false]{\text{\LARGE{$\supseteq$}}}\arrow[draw=none,ul,xshift=0.5cm,"(4)"right]\\
 & \{\text{Serre}\} \arrow[draw=none]{ul}[sloped,auto=false]{\text{\LARGE{$\supseteq$}}}\arrow[draw=none,ul,xshift=-0.5cm,"(5)"left]
 \arrow[draw=none]{ur}[sloped,auto=false]{\text{\LARGE{$\subseteq$}}}\arrow[draw=none,ur,xshift=0.5cm,"(6)"right]& 
	\end{tikzcd}$$
Moreover, each of these inclusions can be proper. Indeed, (2) can be proper by Example~\ref{ex:notTauPerp} and (3) and (4) can be proper by \cite[Example~3.13]{asai_semibricks} and its dual. It is clear that if $\Lambda$ is not $\tau$-tilting finite, then (1) will in general be proper. Likewise, (5) and (6) will be proper unless $\Lambda$ is semisimple or local. We note that all six of these classes coincide when $\Lambda$ is semisimple or local, and that all but the Serre subcategories coincide when $\Lambda$ is $\tau$-tilting finite. Moreover, if $\Lambda$ is hereditary, then (2), (3), and (4) are all equalities by \cite[Corollary~2.17]{IT}.
\end{remark}

%%%%%%%%%%%%%%%%%%%%%%%%%%%%%%%%%%
\section{Reduction of support $\tau$-rigid objects}
\label{sec:reductionFormulas}

The proof of our second main theorem will rely on a reduction formula, which allows us to compare
$\tau$-rigid objects in $\C(\mods\Lambda)$ with those in $\C(\W)$, with $\W\subseteq\mods\Lambda$ a $\tau$-perpendicular subcategory. 
Theorem \ref{thm:BMbijections} below extends Jasso's reduction of $\tau$-rigid modules \cite[Corollary~3.18]{jasso} to account for shifted projectives in the reduction. Note that there is an analogous, and related, formulation in terms of
torsion classes. See \cite[Theorems~3.12 and~3.13]{jasso} and \cite[Theorem~4.12]{DIRRT}. For the statements in this section, we recall the notation for the canonical exact sequence of a module with respect to a torsion pair from Equation~(\ref{eqn:canonical}) in Section \ref{sec:torsion}.

\pagebreak

\begin{theorem}\label{thm:BMbijections}
	Let $\W \subseteq \mods\Lambda$ be a functorially finite wide subcategory.
	\begin{enumerate}
		\item \cite[Proposition~5.6]{BM_exceptional} Let $M \in \W$ be basic and $\tau$-rigid in $\W$. Then there is a bijection
		$$\{V \in \C(\W) \mid M\sqcup V \text{ is basic and $\tau$-rigid in }\C(\W)\}$$
		$$\downarrow\E_M^\W$$
		$$\{V' \in \C(\J_\W(M)) \mid V' \text{ is basic and $\tau$-rigid in }\C(\J_\W(M))\}$$
		summarized as follows.
		\begin{enumerate}
			\item If $N \in \W$ is indecomposable, $M\sqcup N$ is $\tau$-rigid in $\C(\W)$, and $N \notin \Gen M$, then $\E^\W_M(N) = f_{(M^\perp)}(N)$.
			\item If $N \in \W$ is indecomposable, $M\sqcup N$ is $\tau$-rigid in $\C(\W)$, and $N \in \Gen M$, then there exists an indecomposable direct summand $B$ of the Bongartz complement of $M$ (in $\W$) such that $\E^\W_M(N) = f_{(M^\perp)}(B)[1]$.
			\item If $Q \in \P(\W)$ and $M\sqcup Q[1]$ is support $\tau$-rigid in $\C(\W)$, then there exists a direct summand $B$ of the Bongartz complement of $M$ (in $\W$) such that $\E_M^\W(Q[1]) = f_{(M^\perp)}(B)[1]$.
		\end{enumerate}
		The bijection then extends additively.
		\item \cite[Proposition~5.10a]{BM_exceptional} Let $P \in \P(\W)$ be projective in $\W$. Then there is a bijection
		$$\{V \in \C(\W) \mid V\sqcup P[1] \text{ is basic and $\tau$-rigid in }\C(\W)\}$$
		$$\downarrow\E_{P[1]}^\W$$
		$$\{V' \in \C(\J_\W(P[1])) \mid V' \text{ is basic and $\tau$-rigid in }\C(\J_\W(P[1]))\}$$
		given by
		$$\E_{P[1]}^\W(N\sqcup Q[1]) = N\sqcup f_{(P^\perp)}(Q)[1].$$
		\item \cite[Theorem~3.6]{BM_wide} Let $U = M\sqcup P[1] \in \C(\W)$ be basic and support $\tau$-rigid  in $\C(\W)$ and define
		$$\E_{U}^\W:= \E_{\E_M^\W(P[1])}^{\J_\W(M)}\circ \E_M^\W.$$ Then $\E^\W_{M\sqcup P[1]}$ is a bijection
		$$\{V \in \C(\W) \mid U\sqcup V \text{ is basic and support $\tau$-rigid in }\C(\W)\}$$
		$$\downarrow\E_U^\W$$
		$$\{V' \in \C(\J_\W(U)) \mid V' \text{ is basic and support $\tau$-rigid in }\C(\J_\W(U))\}.$$
	\end{enumerate}
\end{theorem}

\begin{remark}\label{rem:BMbijections}
	Since $\E_U^\W$ is defined additively, it follows immediately from Theorem~\ref{thm:rank} that $U \sqcup V$ is support $\tau$-tilting in $\C(\W)$ if and only if $\E_U^\W(V)$ is support $\tau$-tilting in $\J_\W(U)$.
\end{remark}

In case $\W = \mods\Lambda$, we will sometimes denote $\E_U:= \E_U^{\mods\Lambda}$. In Theorem~\ref{thm:associative}, we will show that these ``$\E$-maps'' satisfy an associativity property as follows: If $\W \subseteq\mods\Lambda$ is a $\tau$-perpendicular subcategory of $\mods\Lambda$ and $U\sqcup V$ is basic and support $\tau$-rigid in $\C(\W)$, then $\E_{U\sqcup V}^\W = \E^{\J_\W(U)}_{\E^\W_U(V)}\circ \E^W_U$. This result is established in \cite[Theorem~5.9]{BM_wide} in the case that $\Lambda$ is $\tau$-tilting finite.

For the remainder of this section, we recall some technical results about these bijections and their relationship with Bongartz complements. In essence, these lemmas are well-known (see e.g. \cite[Lemma~4.13]{BM_exceptional}), but we give proofs here for completeness.

\begin{lemma}\label{lem:bongartzGenerator}
	Let $U = M\sqcup P[1] \in \C(\mods\Lambda)$ be support $\tau$-rigid and let $B_U$ be the Bongartz complement of $U$. Then no direct summand of $B_U$ lies in $\Gen M$.
\end{lemma}

\begin{proof}
	Let $B$ be an indecomposable direct summand of $B_U$. By Theorem~\ref{thm:bongartzDef}, we have
	$$\lperp{(\tau M)} \cap P^\perp = \Gen(B_U\sqcup M) = \lperp{(\tau(B\sqcup M))}.$$
	Moreover, by \cite[Lemma~3.8]{BM_wide}, $B\sqcup M$ is $\tau$-rigid in the Serre subcategory $P^\perp$ and satisfies
	$$\lperp{(\tau M)} \cap P^\perp = \lperp{(\tau_{(P^\perp)}M)}\cap P^\perp,\qquad \lperp{(\tau(B\sqcup M))}\cap P^\perp = \lperp{(\tau_{(P^\perp)}(B\sqcup M))}\cap P^\perp.$$
	Since $\Gen(B\sqcup M)\subseteq P^\perp$, it then follows by applying \cite[Proposition~2.22]{AIR} in the category $P^\perp$ that $B \notin \Gen M$.
\end{proof}

\begin{lemma}\label{lem:serreReduction}
	Let $U = M\sqcup P[1] \in \C(\mods\Lambda)$ be  support $\tau$-rigid and basic. Let $B_U$ and $C_U$ be the Bongartz complement and co-Bongartz complement of $U$, respectively. Then
	\begin{enumerate}
		\item For $B \in \add(B_U)$, we have
		$\J(B\sqcup U) = B^\perp \cap \J(U) = (\E_U(B))^\perp \cap \J(U).$
		\item $\P(\J(U)) = \add(\E_{U}(B_U)).$
		\item $\P(\J(U)) = \add(\E_U(C_U)[-1])$.
	\end{enumerate}
\end{lemma}

\begin{proof}
	
	(1) We leave out the proof, since it is identical to the proof of Lemma \ref{lem:serreReduction2}, using Theorem
	\ref{thm:BMbijections}(1)(a) and Lemma \ref{lem:bongartzGenerator}.
	
	(2) By Proposition~\ref{prop:ffReduction}, we note that $\rk(B_U) = \rk(\E_U(B_U)) = \rk(\J(U))$. Thus it suffices to show that $\Gen(\E_U(B_U))\cap \J(U) = \J(U)$. Now let $X \in \J(U)$, so in particular $X \in \lperp{(\tau M)} \cap P^\perp$. By Theorem~\ref{thm:bongartzDef}, we have that $\P(\lperp{(\tau M)})\cap P^\perp = \add(B_U\sqcup M)$. Since $\Hom(M,X) = 0$, this means there exists $B \in \add B_U$ and an epimorphism $g:B\rightarrow X$. Again using the fact that $\Hom(M,X) = 0$, we obtain an epimorphism $g': f_{(M^\perp)}(B) \rightarrow X$. By Theorem~\ref{thm:BMbijections}, we conclude that $f_{(M^\perp)}(B) \in \add(\E_U(B))$ and therefore $X \in \Gen(\E_U(B)) \cap \J(U)$ as claimed.
	
	(3)  By Theorem~\ref{thm:BMbijections}, we have that $\E_U(C_U)$ is support $\tau$-tilting in $\C(\J(U))$ and is contained in $\J(U)[1]$. This means $\E_U(C_U)$ must be the direct sum of the shifts of the indecomposable projectives in $\P(\J(U))$.
\end{proof}

\section{$\tau$-cluster morphism categories}\label{sec:catOfWideSubcats}

In this section, we apply our previous results to extend the definition of the  \emph{$\tau$-cluster morphism category} to arbitrary finite-dimensional algebras. This is a small category whose objects correspond to the $\tau$-perpendicular subcategories of $\mods\Lambda$ and whose morphisms are indexed by support \linebreak $\tau$-rigid pairs in these subcategories. See Definition~\ref{def:tauCMC} below. This category was defined by Igusa and Todorov for hereditary algebras in \cite{IT_exceptional} under the name ``cluster morphism category''. A combinatorial interpretation in 
Dynkin type A using noncrossing partitions and binary forests was also given by Igusa in \cite{igusa_category}. The definition was extended to $\tau$-tilting finite algebras by Marsh and the first author in \cite{BM_wide} under the name ``a category of wide subcategories" and given the name ``$\tau$-cluster morphism category" in \cite{HI_picture}. We state our extension of this definition to arbitrary finite-dimensional algebras now.

\begin{definition}\label{def:tauCMC}
	Let $\Lambda$ be a finite-dimensional algebra. We define the \emph{$\tau$-cluster morphism category} of $\Lambda$, denoted $\mathfrak{W}(\Lambda)$, as follows.
	\begin{enumerate}
		\item The objects of $\mathfrak{W}(\Lambda)$ are the $\tau$-perpendicular subcategories of $\mods\Lambda$.
		\item For $\W\subseteq \mods\Lambda$ a $\tau$-perpendicular subcategory and $U \in \C(\W)$ support  $\tau$-rigid and basic, define a formal symbol $g_U^\W$.
		\item Given $\W_1,\W_2$ two $\tau$-perpendicular subcategories of $\mods\Lambda$, we define
		$$\Hom_{\mathfrak{W}(\Lambda)}(\W_1, \W_2) = \left\{g_U^{\W_1} \middle\vert \begin{array}{l}U \text{ is a basic support $\tau$-rigid object in }\C(\W_1) \\ \text{and }\W_2 = \J_{\W_1}(U)\end{array}\right\}.$$ 
		In particular:
		\begin{enumerate}
			\item If $\W_1 \not\supseteq \W_2$, then $\Hom_{\mathfrak{W}(\Lambda)}(\W_1, \W_2) = \emptyset$.
			\item $\Hom_{\mathfrak{W}(\Lambda)}(\W_1, \W_1) = g_0^{\W_1}$
		\end{enumerate}
		\item Given $g_U^{\W_1}:\W_1 \rightarrow \W_2$ and $g_V^{\W_2}:\W_2\rightarrow \W_3$ in $\mathfrak{W}(\Lambda)$, denote $\widetilde{V}:= \left(\E_U^{\W_1}\right)^{-1}(V)$. We define
		$$g_V^{\W_2} \circ g_U^{\W_1} = g_{U\sqcup \widetilde{V}}^{\W_1}.$$
	\end{enumerate}
\end{definition}

\begin{remark}
For $\tau$-tilting finite algebras, it is well-known 
\cite{MS}
that all wide subcategories are both left and right finite, so in particular they are $\tau$-perpendicular. The above definition therefore
specializes to the definition in~\cite{BM_wide}.
\end{remark}

\begin{remark}
    An independent generalization of the $\tau$-cluster morphism category to arbitrary finite-dimensional algebras is given in the concurrent work of B{\o}rve \cite{borve}. The construction given there replaces $\tau$-perpendicular subcategories with certain thick subcategories of the bounded derived category $\Db(\mods\Lambda)$ and replaces support $\tau$-rigid objects with 2-term presilting objects. The composition law can then be described in terms of the (pre)silting reduction of Iyama--Yang \cite{IY}. It is shown explicitly in \cite{borve} that our generalization and B{\o}rve's yield categories which are equivalent.
\end{remark}

The main goal of this section is to prove that the $\tau$-cluster morphism category is indeed a well-defined category (Theorem~\ref{thm:intro:mainB}, restated as Theorem~\ref{thm:mainB} below). As with \cite{IT_exceptional,igusa_category,BM_wide}, the main technicality is in showing that the composition law is well-defined and associative.
In the present paper, this will be a consequence of the following generalization of \cite[Theorem~4.3]{BM_wide}. 

\begin{theorem}\label{thm:composition}
	Let $\W\subseteq\mods\Lambda$ be a $\tau$-perpendicular subcategory of $\mods\Lambda$. Let $U\sqcup V$ be basic and support $\tau$-rigid in $\C(\W)$. Then
	$$\J_\W(U\sqcup V) = \J_{\J_{\W}(U)}(\E^\W_U(V)).$$
\end{theorem}

Our proof of Theorem~\ref{thm:composition} is largely contained in the two technical lemmas which follow.

\begin{lemma}\label{lem:associative1}
	Let $U \in \C(\mods\Lambda)$ and $N\in\mods\Lambda$ such that $U\sqcup N$ is support $\tau$-rigid and basic. Write $U = M \sqcup P[1]$ and let $\overline{N}$ be the direct sum of the indecomposable direct summands of $N$ which do not lie in $\Gen M$. Then the following coincide:
	\begin{enumerate}
		\item $(\Gen N)\cap \J(U)$
		\item $(\Gen(M\sqcup N))\cap \J(U)$
		\item $f_{(M^\perp)}(\Gen(M\sqcup N))$
		\item $(\Gen(f_{(M^\perp)}(N))) \cap \J(U)$
		\item $(\Gen(\E_U(\overline{N})))\cap \J(U)$
	\end{enumerate}
\end{lemma}

\begin{proof}
	The equality $(1) = (2)$ follows immediately from the fact that $\J(U) \subseteq M^\perp$.

	We next show that $(2) = (3)$. Note that by definition $$f_{(M^\perp)}(\Gen(M\sqcup N)) \subseteq (\Gen(M\sqcup N))\cap M^\perp.$$ Moreover, we have that $\Gen(M\sqcup N) \subseteq \lperp{(\tau M)}\cap P^\perp$ since $U\sqcup N$ is support   $\tau$-rigid. Now, if $X \in (\Gen(M\sqcup N))\cap \J(U)$, then in particular $X \in M^\perp$ and so $f_{(M^\perp)}(X) = X$. We conclude that $(2) = (3)$.

	We now show that $(3) = (4)$. It is shown in \cite[Lemma~5.5]{BM_wide} that
	$$f_{(M^\perp)}(\Gen(M\sqcup N)) = (\Gen(f_{(M^\perp)}(N))) \cap M^\perp \cap \lperp{(\tau M)}.$$
	Moreover, since $N \in P^\perp$ and $f_{(M^\perp)}(N) \in \Gen N$, we have that $\Gen(f_{(M^\perp)}(N)) \subseteq P^\perp$. It follows that $(3) = (4)$.
	
	It remains to show that $(4) = (5)$. This follows from the definition of $\E_U$ (see Theorem~\ref{thm:BMbijections}) and the fact that $f_{(M^\perp)}(N) = f_{(M^\perp)}(\overline{N})$.
\end{proof}

\begin{lemma}\label{lem:associative2}
	Let $U \sqcup V \in \C(\mods\Lambda)$ be support $\tau$-rigid and basic. Let $B$ be the Bongartz complement of $U\sqcup V$ (in $\mods\Lambda$). Write $U = M\sqcup P[1], V = N\sqcup Q[1]$, and $\E_U(V) = L\sqcup R[1]$. Then the following coincide:
	\begin{enumerate}
		\item $\lperp{(\tau N)}\cap Q^\perp \cap \J(U)$
		\item $\lperp{(\tau N\sqcup \tau M)} \cap (Q\sqcup P)^\perp \cap \J(U)$
		\item $\Gen(B\sqcup N \sqcup M)\cap \J(U)$
		\item $\Gen(\E_U(B)\sqcup L)\cap \J(U)$
		\item $\lperp{(\tau_{\J(U)}L)}\cap R^\perp \cap \J(U)$
	\end{enumerate}
\end{lemma}

\begin{proof}
	The equality $(1) = (2)$ follows immediately from the fact that $\J(U) \subseteq \lperp{(\tau M)} \cap P^\perp$. Likewise, the equality $(2) = (3)$ follows immediately from the definition of the Bongartz complement (Theorem~\ref{thm:bongartzDef}).

We next show that $(3) = (4)$.  By Theorem~\ref{thm:BMbijections} and Lemma~\ref{lem:bongartzGenerator}, we note that $\E_U(B) \in \mods\Lambda$. Moreover, Theorem~\ref{thm:BMbijections} also implies (see Remark~\ref{rem:BMbijections}) that $\E_U(B) \sqcup L \sqcup R[1] =
\E_U(B \sqcup V)$
is support $\tau$-tilting in $\C(\J(U))$. 
This means $(3) = (4)$ is a special case of equation $(2) = (5)$ in Lemma~\ref{lem:associative1}. 

We proceed to show that $(4) \subseteq (5)$.
First note that since $R$ is projective in $\J(U)$, we have that $\lperp{(\tau_{\J(U)}L)}\cap R^\perp \cap \J(U)$ is closed under factors in $\J(U)$, and hence it suffices to show that 
$\E_U(B)\sqcup L \in \lperp{(\tau_{\J(U)}L)}\cap R^\perp \cap \J(U)$. This follows from the fact that  
$\E_U(B) \sqcup L \sqcup R[1]$ is support $\tau$-tilting in $\C(\J(U))$.
	
	We will conclude by showing that $(5) \subseteq (2)$. Let $B'$ be the Bongartz complement of $\E_U(V)$ in $\J(U)$. Since $B'$ is a module, Theorem~\ref{thm:BMbijections} implies that $\overline{B}:= \E_U^{-1}(B') \in \mods\Lambda$ and that $\overline{B}\sqcup U \sqcup V$ is support $\tau$-tilting. 
	We claim that 
	\begin{eqnarray*}
		\lperp{(\tau_{\J(U)}L)}\cap R^\perp \cap \J(U) &=& \Gen(B'\sqcup L) \cap \J(U)\\
			&=& \Gen(\overline{B}\sqcup N \sqcup M) \cap \J(U)\\
			&\subseteq& \lperp{(\tau N\sqcup \tau M)} \cap (Q\sqcup P)^\perp \cap \J(U)
	\end{eqnarray*}
	The first equality follows from Theorem~\ref{thm:bongartzDef},
	and the second from the equality $(2)= (5)$ in
	Lemma~\ref{lem:associative1}. The inclusion follows from the 
	fact that $\overline{B}\sqcup U \sqcup V$ is support $\tau$-tilting, using that $(2)$ is closed under factors in 
	$\J(U)$.
\end{proof}

We now proceed with our proof of Theorem~\ref{thm:composition}

\begin{proof}[Proof of Theorem~\ref{thm:composition}]
	Since $\W$ is equivalent to the module category of a basic finite-dimensional algebra, it suffices to consider the case where $\W = \mods\Lambda$. Write $U = M\sqcup P[1]$, $V = N\sqcup Q[1]$, and $\E_U(V) = L\sqcup R[1]$. For readability, denote $\T_0 = \Gen(N\sqcup M)$ and $\T_1 = \lperp{(\tau N\sqcup \tau M)}\cap (Q\sqcup P)^\perp$. Noting that $\J(U\sqcup V) \subseteq \J(U)$, 
	Lemmas~\ref{lem:associative1} and~\ref{lem:associative2} then imply that
	\begin{eqnarray*}
		\J(U\sqcup V) &=& (\T_0^\perp \cap \J(U)) \cap (\T_1 \cap \J(U))\\
			&=& (\Gen L)^\perp \cap \lperp{(\tau_{\J(U)} L)}\cap R^\perp \cap \J(U)\\
			&=& L^\perp \cap \lperp{(\tau_{\J(U)} L)}\cap R^\perp \cap \J(U)\\
			&=& \J_{\J(U)}(\E_U(V)).
	\end{eqnarray*}
\end{proof}

Before we proceed with proving the second main theorem, we note that Theorem \ref{thm:composition} has some 
interesting consequences.

\begin{corollary}[Corollary~\ref{cor:intro:main}]\label{cor:mainB}
	Let $\Lambda$ be a finite-dimensional algebra. Let $\V\subseteq \W \subseteq\mods\Lambda$ be a chain of subcategories such that $\V$ is a $\tau$-perpendicular subcategory of $\W$ and $\W$ is a $\tau$-perpendicular subcategory of $\mods\Lambda$. Then $\V$ is a $\tau$-perpendicular subcategory of $\mods\Lambda$.
\end{corollary}

\begin{proof}
	Let $U \in \C(\mods\Lambda)$ and $V \in \C(\W)$ such that $\W = \J(U)$ and $\V = \J_\W(V)$. By Theorem~\ref{thm:BMbijections} and Theorem~\ref{thm:composition}, it follows that $U \sqcup \left(\E_U^{-1}(V)\right) \in \C(\mods\Lambda)$ is support $\tau$-rigid and satisfies
	$$\J\left(U \sqcup \left(\E_U^{-1}(V)\right)\right) = \J_{\J(U)}\left(\E_U\circ\E_U^{-1}(V)\right) = \J_{\W}(V) = \V.$$
	We conclude that $\V$ is a $\tau$-perpendicular subcategory of $\mods\Lambda$.
\end{proof}

In many cases, the converse of Corollary~\ref{cor:intro:main} holds as well. For example, in the $\tau$-tilting finite case all wide subcategories are $\tau$-perpendicular and in the hereditary case, $\tau$-perpendicular subcategories and left finite wide subcategories coincide. Each of these implies that if $\V$ and $\W$ are $\tau$-perpendicular subcategories of $\mods\Lambda$ with $\V \subseteq \W$, then $\V$ is a $\tau$-perpendicular subcategory of $\W$. We expect that this is the case in general; that is, we propose the following conjecture.

\begin{conjecture}\label{conj:iteratedReduction}
    Let $\Lambda$ be a finite-dimensional algebra. Let $\W \subseteq\mods\Lambda$ be a $\tau$-perpendicular subcategory of $\mods\Lambda$ and let $\V \subseteq \W$ be a wide subcategory of $\W$. Then $\V$ is a $\tau$-perpendicular subcategory of $\mods\Lambda$ if and only if $\V$ is a $\tau$-perpendicular subcategory of $\W$.
\end{conjecture}

As another consequence of Theorem~\ref{thm:composition}, we have the following.

\begin{corollary}\label{cor:iteratedLeftFinite}
	Let $\Lambda$ be a finite-dimensional algebra and let $\W \subseteq \mods\Lambda$ be a subcategory. Then $\W$ is a $\tau$-perpendicular subcategory of $\mods\Lambda$ if and only if there exists a subcategory $\V$ with $\W\subseteq \V\subseteq \mods\Lambda$ such that $\V$ is a left finite wide subcategory of $\mods\Lambda$ and $\W$ is a left finite wide subcategory of $\V$. Moreover, the statement is true if one or both instances of ``left'' are replaced with ``right''.
\end{corollary}

\begin{proof}
	First suppose $\W = \J(U)$ is a $\tau$-perpendicular subcategory of $\mods\Lambda$. By Theorem~\ref{thm:mainA}, there exists $\W \subseteq \V \subseteq \mods\Lambda$ such that $\V$ is a left finite wide subcategory of $\mods\Lambda$ and $\W$ is a Serre subcategory of $\V$. By Corollary~\ref{cor:serreFinite}, it follows that $\W$ is a left finite wide subcategory of $\V$ as well.
	
	Now suppose that there exists $\W \subseteq \V \subseteq \mods\Lambda$ such that $\W$ is a left finite wide subcategory of $\V$ and $\V$ is a left finite wide subcategory of $\mods\Lambda$. Then $\W$ is a $\tau$-perpendicular subcategory of $\mods\Lambda$ as an immediate consequence of Lemma~\ref{lem:finiteIsJasso} and Corollary~\ref{cor:mainB}.
	
	The proofs where one or both instances of ``left'' are replaced with ``right'' are identical.
\end{proof}

We now proceed with the proof of the main theorem. In \cite{BM_wide}, which deals with the 
$\tau$-tilting finite case,
one establishes associativity of the composition operation, by proving that 
\begin{equation}\label{eq:assos} 
\E_{\E_U(V)}^{\J(U)} \circ \E_U = \E_{U\sqcup V}   
\end{equation}
for any basic support $\tau$-rigid $U\sqcup V \in \C(\mods\Lambda)$.
This is shown to be a consequence of the fact that 
\begin{equation}\label{eq:comp} 
\J_\W(U\sqcup V) = \J_{\J_{\W}(U)}(\E^\W_U(V))
\end{equation}
However, the proof that (\ref{eq:comp}) implies (\ref{eq:assos}) given in \cite[Sections 5-9]{BM_wide} does
not use that $\Lambda$ is $\tau$-tilting finite. We have shown in Theorem \ref{thm:composition} that 
Equation (\ref{eq:comp}) holds for $\tau$-perpendicular subcategories in the general case, and
hence we obtain (\ref{eq:assos}) for free. Note that in the $\tau$-tilting finite case, in fact all
wide subcategories are $\tau$-perpendicular.

We will here provide an alternative and much more efficient proof of why (\ref{eq:comp}) implies (\ref{eq:assos}) in the general case, 
which only builds on two short lemmas in \cite{BM_wide}, namely Lemmas~5.5 (via Lemma~\ref{lem:associative1} in the present paper) and~6.2.
This is completed in Theorem \ref{thm:associative} below, but we first prepare with an additional technical lemma.

\begin{lemma}\label{lem:associative3}
	Let $U\sqcup V \in \C(\mods\Lambda)$ and $L \in \mods\Lambda$ such that $U\sqcup V \sqcup L$ is support $\tau$-rigid and basic. Let $\overline{N}$ be the direct sum of the indecomposable direct summands of $N$ which do not lie in $\Gen M$. Let $L'$ be an indecomposable direct summand of $L$. Then the following are equivalent.
	\begin{enumerate}
		\item $\E_{U\sqcup V}(L')$ is a module.
		\item $L' \notin \Gen(M\sqcup N)$
		\item $L' \notin \Gen M$ and $f_{(M^\perp)}(L') \notin (\Gen(f_{(M^\perp)}(\overline{N})))\cap \J(U)$
		\item $\E_U(L')$ is a module and $\E_U(L') \notin (\Gen(\E_U(\overline{N})))\cap \J(U)$
		\item $\E_{\E_U(V)}^{\J(U)} \circ \E_{U}(L')$ is a module.
	\end{enumerate}
Moreover, if (1)-(5) hold then $\E_{\E_U(V)}^{\J(U)} \circ \E_{U}(L') = \E_{U\sqcup V}(L')$.
\end{lemma}

\begin{proof}
	The equivalences $(1\iff 2)$, $(3\iff 4)$, and $(4\iff 5)$ all follow from the definitions of the ``$\E$-maps'' given in Theorem~\ref{thm:BMbijections}.
	
	For the equivalence $(2\iff 3)$, we note that $L' \notin \Gen(M\sqcup N)$ if and only if
	$$f_{(M^\perp)}(L') \notin f_{(M^\perp)}(\Gen(M\sqcup N)) = (\Gen(f_{(M^\perp)}(\overline{N}))\cap \J(U)$$
	by Lemma~\ref{lem:associative1}.
	This, together with the fact that if $L' \notin \Gen(M\sqcup N)$ then $L' \notin \Gen M$, proves the equivalence of (2) and (3).
	
	Now suppose that (1)-(5) hold. For readability, denote $L'':= f_{(M^\perp)}(L')$ and denote $\T := f_{(M^\perp)}(\Gen(M\sqcup N))$. Recall from Lemma~\ref{lem:associative1} that $\T = (\Gen(f_{(M^\perp)}(N))) \cap \J(U)$. In particular, we have $\T \subseteq \J(U)\subseteq \J(M)$. Now denote by
	$$0\rightarrow t_\T(L'') \rightarrow L'' \rightarrow f_{(\T^\perp\cap \J(U))}(L'')\rightarrow 0$$
	the canonical exact sequence with respect to the torsion pair $(\T,\T^\perp \cap \J(U))$ in $\J(U)$. Likewise, denote by
	$$0\rightarrow t_\T(L'') \rightarrow L'' \rightarrow f_{(\T^\perp\cap \J(M))}(L'')\rightarrow 0$$
	the canonical exact sequence with respect to the torsion pair $(\T,\T^\perp \cap \J(M))$ in $\J(M)$. Since both sequences start with $t_\T(L'')$, we see that $f_{(\T^\perp\cap \J(U))}(L'') = f_{(\T^\perp\cap \J(M))}(L'')$.
	
	We now observe that $\E_{\E_U(V)}^{\J(U)} \circ \E_{U}(L') = f_{(\T^\perp\cap\J(U))}(L'')$ and $\E_{(U\sqcup V)} = f_{((M\sqcup N)^\perp)}(L')$ by construction. Finally, it is shown in \cite[Lemma~6.2]{BM_wide} that $f_{(\T^\perp\cap \J(M))}(L'') = f_{((M\sqcup N)^\perp)}(L')$. We conclude that $\E_{\E_U(V)}^{\J(U)} \circ\E_{U}(L') = \E_{U\sqcup V}(L')$ as desired.
\end{proof}

\begin{remark}\label{rem:associative3}
	Since the bijections $\E_U$ and $\E^\W_U$ are additive, the assumption that $L'$ is indecomposable in Lemma~\ref{lem:associative2} can be replaced with the assumption that no direct summand of $L'$ lies in $\Gen(M\sqcup N)$.
\end{remark}

We are now prepared to verify Equation~(\ref{eq:assos}).

\begin{theorem}\label{thm:associative}
	Let $\W\subseteq\mods\Lambda$ be a $\tau$-perpendicular subcategory of $\mods\Lambda$ and let $U\sqcup V \in \C(\W)$ be support $\tau$-rigid and basic. Then $$\E_{\E_U(V)}^{\J_\W(U)} \circ \E^\W_U = \E^\W_{U\sqcup V}.$$
\end{theorem}

\begin{proof}
	Since $\W$ is equivalent to the module category of a basic finite-dimensional algebra, it suffices to consider the case where $\W = \mods\Lambda$.

	Let $W \in \C(\mods\Lambda)$ such that $U\sqcup V\sqcup W$ is support $\tau$-rigid and basic. Write $U = M\sqcup P[1]$ and $V = N\sqcup Q[1]$. Let $\overline{L}$ be the direct sum of the indecomposable direct summands of $W$ which are modules and do not lie in $\Gen(M\sqcup N)$, and let $W' \in \C(\mods\Lambda)$ such that $\overline{L}\sqcup W' = W$. Let $B \in \mods\Lambda$ be the Bongartz complement of $U\sqcup V\sqcup W$. We recall from Lemma~\ref{lem:bongartzGenerator} that no direct summand of $B$ lies in $\Gen(M\sqcup N)$. Therefore, by Lemma~\ref{lem:associative3} and Remark~\ref{rem:associative3}, we have that $\E_{\E_U(V)}^{\J(U)} \circ \E_U(B\sqcup\overline{L}) = \E_{U\sqcup V}(B\sqcup\overline{L})$. Moreover, this equation also holds if $B\sqcup \overline{L}$ is replaced with any of its (not necessarily indecomposable) direct summands.
	
	For readability, denote $B' := \E_{U\sqcup V}(B)$ and $L' := \E_{U\sqcup V}(\overline{L})$. Now recall from Theorem~\ref{thm:composition} that $\J(U\sqcup V) = \J_{\J(U)}(\E_U(V))$. Theorem~\ref{thm:bongartzDef} and Lemma~\ref{lem:associative2} then imply that (i) both $\E_{U\sqcup V}(W')$ and $\E_{\E_U(V)}^{\J(U)} \circ \E_U(W')$ lie in $\J(U\sqcup V)[1]$, and (ii) both $B'\sqcup L' \sqcup \E_{U\sqcup V}(W')$ and $B'\sqcup L'\sqcup \E_{\E_U(V)}^{\J(U)} \circ \E_U(W')$ are support $\tau$-tilting in $\C(\J(U\sqcup V))$. This implies that $\E_U(W') = \E_{\E_U(V)}^{\J(U)} \circ \E_U(W')$, as both coincide with the co-Bongartz complement of $B'\sqcup L'$ in $\J(U\sqcup V)$.
	
Now recall that $W = \overline{L}\sqcup W'$. Since all of the bijections $\E_U$, $\E_{U\sqcup V}$, and $\E_{\E_U(V)}^{\J(U)}$ are additive, the previous two paragraphs imply that $\E_{\E_U(V)}^{\J(U)}\circ \E_U(W) = \E_{U\sqcup V}(W)$. This completes the proof.
\end{proof}

We are now ready to complete the proof that the $\tau$-cluster morphism category is indeed a category. This essentially follows from Theorem~\ref{thm:associative} identically as in \cite[Corollary~1.10]{IT_exceptional}, \cite[Section~1]{igusa_category}, and \cite[Corollary~1.8]{BM_wide}.

\begin{theorem}[Theorem~\ref{thm:intro:mainB}]\label{thm:mainB}
	Let $\Lambda$ be a finite-dimensional algebra. Then $\mathfrak{W}(\Lambda)$ is a well-defined category.    
\end{theorem}

\begin{proof}
	It is straightforward to show that for any $\tau$-perpendicular subcategory $\W \subseteq \mods\Lambda$, the morphism $g_0^\W$ is the identity of $\W$. Thus we need only show that the composition law is associative. Consider
	$$\W_1\xrightarrow{g_U^{\W_1}} \W_2 \xrightarrow{g_V^{\W_2}} \W_3 \xrightarrow{g_W^{\W_3}} \W_4$$
a sequence of three composable morphisms in $\mathfrak{W}$. Then, by Theorem~\ref{thm:associative} and the additivity of the ``$\E$-maps'', we have:
\begin{eqnarray*}
	\left(g_W^{\W_3}\circ g_{V}^{\W_2}\right)\circ g_U^{\W_1} &=& g_{V\sqcup \left(\E_V^{\W_2}\right)^{-1}(W)}^{\W_2}\circ g_U^{\W_1}\\
	&=& g_{U\sqcup \left(\E_U^{\W_1}\right)^{-1}(V) \sqcup \left(\E_V^{\W_2}\circ \E_U^{\W_1}\right)^{-1}(W)}^{\W_1}\\
	&=& g_W^{\W_3} \circ g_{U\sqcup \left(\E_W^{\W_1}\right)^{-1}(V)}^{\W_1}\\
	&=& g_W^{\W_3}\circ\left( g_{V}^{\W_2}\circ g_W^{\W_1}\right)
\end{eqnarray*}
\end{proof}

For $\W\subseteq \mods\Lambda$ a $\tau$-perpendicular subcategory, we can likewise define the  $\tau$-cluster morphism category $\mathfrak{W}(\W)$ in the usual way; i.e., by identifying $\W$ with some module category. We then obtain the following.

\begin{proposition}\label{prop:reductionSubcat}
	Let $\W \subseteq \mods\Lambda$ be a $\tau$-perpendicular subcategory. Then $\mathfrak{W}(\W)$ is equivalent to the full subcategory of $\mathfrak{W}(\Lambda)$ whose objects are the $\tau$-perpendicular subcategories of $\W$.
\end{proposition}

\begin{proof}
	This is a straightforward consequence of Corollary~\ref{cor:mainB}.
\end{proof}

We conclude this section by generalizing the results of \cite[Section~10]{BM_wide}.

\begin{proposition}\label{prop:distinctReductions2}\
	\begin{enumerate}
		\item Let $M$ be an indecomposable non-projective $\tau$-rigid module and let $B_M$ be the Bongartz complement of $M$. Then $M \in \Gen B_M$ and $\J(M) = \W_L(\Gen B_M)$.
		\item Let $M$ and $N$ be indecomposable $\tau$-rigid modules. Then $\J(M) = \J(N)$ if and only if $M \cong N$.
	\end{enumerate}
\end{proposition}

\begin{proof}
	(1) As in the proof of \cite[Lemma 10.6]{BM_wide}, we have that the indecomposable direct summands of $B$ are split projective in $\lperp(\tau M)$. If in addition $M$ is split projective in $\lperp(\tau M)$, then $\Gen(B \sqcup M) = \mods\Lambda$ by Theorem~\ref{thm:bongartzDef} and Proposition~\ref{prop:bongartzSplitNonsplit}. Since $B\sqcup M$ is $\tau$-tilting, this implies that $B\sqcup M = \Lambda$ and $M$ is projective, a contradiction. We conclude that $M$ is not split projective in $\lperp(\tau M)$, and so the result is a special case of Lemma~\ref{lem:finiteIsJasso}.
	
	(2) Suppose that $\J(M) = \J(N)$. We first consider the case where neither $M$ nor $N$ is projective. Let $B_M$ be the Bongartz complement of $M$ and $B_N$ the Bongartz complement of $N$. Then by (1) we have $\W_L(\Gen B_M) = \J(M) = \J(N) = \W_L(\Gen B_N)$. This then implies that $\Gen B_M = \Gen B_N$. It follows that $M \cong N$ is the unique indecomposable ext-projective in this torsion class which is not split-projective.
	
	Now suppose that $M$ is projective. Since $M$ is indecomposable, we note that $\mathsf{top}(M)$ is simple. Moreover, given an arbitrary simple $S \in \mods\Lambda$, we have that $S \in \J(M)$ if and only if $S \ncong \mathsf{top}(M)$. In particular, if $S \ncong \mathsf{top}(M)$, then $\Hom(N,S) = 0 = \Hom(S,\tau N)$. We conclude that $\mathsf{top}(M) \cong \mathsf{top}(N)$ and that $\mathsf{soc}(\tau N) = 0$. In particular, this means $M \cong N$.
\end{proof}

\begin{theorem}
    Let $\Lambda$ be a finite-dimensional algebra and let $\W \subseteq \mods\Lambda$ be a $\tau$-perpendicular subcategory. Let $\V \subseteq \W$ be a $\tau$-perpendicular subcategory of $\W$ such that $\rk(\W) = \rk(\V) + 1$. Then exactly one of the following occurs:
    \begin{enumerate}
        \item There is exactly one morphism in $\mathfrak{W}(\Lambda)$ from $\W$ to $\V$ and $\V = \J_\W(M)$ for some indecomposable module $M$ which is $\tau$-rigid, but not projective, in $\W$.
        \item There are exactly two morphisms in $\mathfrak{W}(\Lambda)$ from $\W$ to $\V$ and $\V = \J_\W(P) = \J_\W(P[1])$ for some indecomposable module $P$ which is projective in $\W$.
    \end{enumerate}
\end{theorem}
\begin{proof}
    By Proposition~\ref{prop:reductionSubcat}, we can assume that $\W = \mods\Lambda$. Theorem~\ref{thm:rank} then implies that there exists an indecomposable support $\tau$-rigid $U \in \C(\mods\Lambda)$ such that $\J(U) = \V$. The result then follows from Proposition~\ref{prop:distinctReductions2} and that fact that $\J(P) = \J(P[1])$ for any projective module $P$.
\end{proof}

\begin{remark}\label{rem:notLattice?}
        In \cite[Section~10]{BM_wide}, the results generalized here are stated as relationships between the lattice of wide subcategories and the $\tau$-cluster morphism category. (The partial order on the lattice of wide subcategories is given by containment,  and the meet of two wide subcategories is their intersection.) In the $\tau$-tilting infinite case, however, the set of $\tau$-perpendicular subcategories may not be closed under intersections and may not form a lattice. See for example \cite[Examples~3.2.2 and~3.2.3]{ringel}, which show that path algebras of type $\widetilde{A}_{2,1}$ and $\widetilde{A}_{2,2}$, respectively, exhibit such behavior\footnote{In these examples, the term ``exceptional subcategory'' is used to mean functorially finite wide subcategory. We also recall from Remark~\ref{rem:inclusions} that over hereditary algebras, these are precisely the $\tau$-perpendicular subcategories.}. The authors thank Haruhisa Enomoto for sharing these examples with them.
\end{remark}

%%%%%%%%%%%%%%%%%%%%%%%%%%%%%%%%

\section{An example}

In this section, we consider a pair of examples. As examples in the $\tau$-tilting finite case can be found in \cite[Section~12]{BM_wide}, our examples are both $\tau$-tilting infinite. This means the corresponding $\tau$-cluster morphism categories are infinite as well.

We first consider the Kronecker path algebra $\Lambda_1 = K(1\rightrightarrows 2).$
   For $i \in \mathbb{N}$, we denote by $M_{(i,i+1)}$ and $M_{(i+1,i)}$ the unique (up to isomorphism) indecomposable $\Lambda_1$-modules with dimension vectors $(i,i+1)$ and $(i+1,i)$. We note that $P(1) = M_{(1,2)}$, $S(1) = M_{(1,0)}$, and $P(2) = S(2) = M_{(0,1)}$.
    
    Figure~\ref{fig:kronecker} is an illustration of the category $\mathfrak{W}(\Lambda_1)$. The vertices are the  $\tau$-perpendicular subcategories of $\mods\Lambda$ (which in this case are precisely the functorially finite wide subcategories). An irreducible morphism $g_U^\W:\W\rightarrow \V$ (so that $U$ is indecomposable and support $\tau$-rigid in $\C(\W)$ and $\J_\W(U) = \V)$ is shown as an arrow $\W\rightarrow \V$ labeled by $U$. The wide subcategories $\add\{M_{(i,i+1)}\}$ generated by the preprojective modules all appear above the horizontal dashed line, with $i$ increasing as one moves counter clockwise. Likewise, the wide subcategories $\add\{M_{(i+1,i)}\}$ generated by the preinjective modules all appear below the horizontal dashed line, with $i$ increasing as one moves clockwise. The category is drawn so that every square commutes, and wide subcategories which appear more than once in the figure should be identified.

\begin{center}
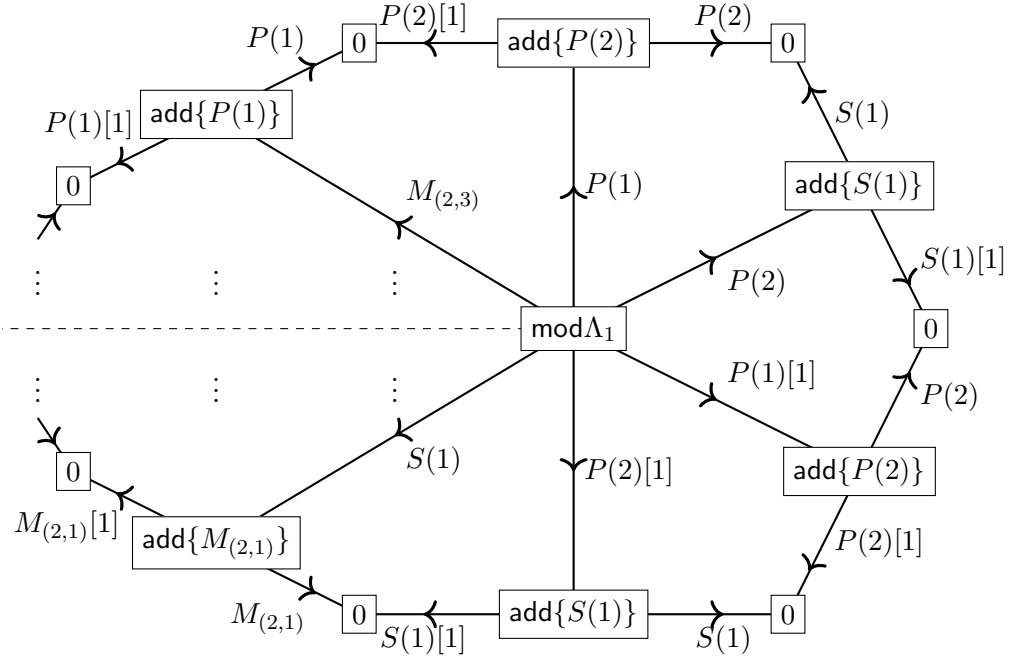
\begin{figure}
\begin{tikzpicture}[scale=0.95]
\begin{scope}[thick,decoration={
	markings,
	mark=at position 0.5 with {\arrow[scale=1.5]{>}}}
	]
        \draw[postaction={decorate}] (0,0)--(0,4) node [midway,anchor=west]{$P(1)$};
        \draw[postaction={decorate}](0,0)--(-5,3) node [midway,anchor=south west]{$M_{(2,3)}$};
        \draw[postaction={decorate}](-7.5,1.25)--(-7,2);
        \draw[postaction={decorate}] (0,0)--(0,-4) node [midway,anchor=west]{$P(2)[1]$};
        \draw[postaction={decorate}](0,0)--(-5,-3) node [midway,anchor=north west]{$S(1)$};
        \draw[postaction={decorate}](-7.5,-1.25)--(-7,-2);
        \draw[postaction={decorate}] (0,0)--(4,2) node [midway,anchor=north west]{$P(2)$};
        \draw[postaction={decorate}] (0,0)--(4,-2) node [midway,anchor=south west]{$P(1)[1]$};
    \end{scope}
	
\begin{scope}[thick,decoration={
	markings,
	mark=at position 0.7 with {\arrow[scale=1.5]{>}}}
	]
        \draw[postaction={decorate}] (0,4)--(-3,4)node [pos=0.7,anchor=south]{$P(2)[1]$};
        \draw[postaction={decorate}](-5,3)--(-3,4)node [pos=0.7,anchor=south east]{$P(1)$};
        \draw[postaction={decorate}](-5,3)--(-7,2)node [midway,anchor=south east]{$P(1)[1]$};
        \draw[postaction={decorate}] (0,-4)--(-3,-4)node [pos=0.7,anchor=north]{$S(1)[1]$};
        \draw[postaction={decorate}](-5,-3)--(-3,-4)node [pos=0.7,anchor=north east]{$M_{(2,1)}$};
        \draw[postaction={decorate}](-5,-3)--(-7,-2)node [pos=0.6,anchor=north east]{$M_{(2,1)}[1]$};
        \draw[postaction={decorate}](4,2)--(3,4)node [midway,anchor= west]{$S(1)$};
        \draw[postaction={decorate}](4,-2)--(3,-4)node [midway,anchor= west]{$P(2)[1]$};
        \draw[postaction={decorate}](0,-4)--(3,-4)node [pos=0.7,anchor= north]{$S(1)$};
        \draw[postaction={decorate}](0,4)--(3,4)node [pos=0.7,anchor= south]{$P(2)$};
        \draw[postaction={decorate}](4,2)--(5,0)node [pos=0.7,anchor= south west]{$S(1)[1]$};
        \draw[postaction={decorate}](4,-2)--(5,0)node [pos=0.7,anchor= north west]{$P(2)$};
\end{scope}
        \draw[dashed] (0,0)--(-8,0);
        
\node at (0,0) [draw,fill=white] {$\mods\Lambda_1$};
    
\node at (-3,4) [draw,fill=white] {$0$};
\node at (-7,2) [draw,fill=white] {$0$};
\node at (-3,-4) [draw,fill=white] {$0$};
\node at (-7,-2) [draw,fill=white] {$0$};
\node at (3,4) [draw,fill=white] {$0$};
\node at (3,-4) [draw,fill=white] {$0$};
\node at (5,0) [draw,fill=white] {$0$};

\node at (0,4) [draw,fill=white] {$\add\{P(2)\}$};
\node at (-5,3) [draw,fill=white] {$\add\{P(1)\}$};
\node at (0,-4) [draw,fill=white] {$\add\{S(1)\}$};
\node at (-5,-3) [draw,fill=white] {$\add\{M_{(2,1)}\}$};
\node at (4,2) [draw,fill=white] {$\add\{S(1)\}$};
\node at (4,-2) [draw,fill=white] {$\add\{P(2)\}$};

\node at (-7.5,0.75){$\vdots$};
\node at (-5,0.75){$\vdots$};
\node at (-2.5,0.75){$\vdots$};
\node at (-7.5,-0.75){$\vdots$};
\node at (-5,-0.75){$\vdots$};
\node at (-2.5,-0.75){$\vdots$};

\end{tikzpicture}
\caption{The category $\mathfrak{W}(\Lambda_1)$ for $\Lambda_1 = K(1\rightrightarrows 2)$.}\label{fig:kronecker}
\end{figure}
\end{center}

%%%%%%%%%%%%%%%%%%%%%%%%%%%

For our second example, we consider quiver $Q = 1\rightrightarrows 2 \rightarrow 3$ and the algebra $\Lambda_2 = KQ/\mathrm{rad}^2 KQ$. Again for $i \in \mathbb{N}$, we denote by $M_{(i,i+1,0)}$ and $M_{(i+1,i,0)}$ the unique (up to isomorphism) $\Lambda_2$-modules with dimension vectors $(i,i+1,0)$ and $(i+1,i,0)$. The irreducible morphisms in $\mathfrak{W}(\Lambda_2)$ with source $\mods\Lambda_2$ are shown in Figure~\ref{fig:example}. Similarly to before, a morphism $g_U^{\mods\Lambda_2}: \mods\Lambda_2 \rightarrow \W$ is labeled by $U$. Moreover, every module of the form $M_{(i,i+1,0)}$ or $M_{(i+1,i,0)}$ corresponds to some morphism with source $\mods\Lambda_2$.

To complete the picture, we can utilize Proposition~\ref{prop:reductionSubcat}. The Serre subcategory $P(3)^\perp$ is equivalent to $\mods\Lambda_1$, so there is a copy of $\mathfrak{W}(\Lambda_1)$ sitting inside of $\mathfrak{W}(\Lambda_2)$ which has $P(3)^\perp$ identified with $\mods\Lambda_1$. The Serre subcategory $P(1)^\perp = \add\{P(2),S(2),P(3)\}$ is equivalent to the module category of the path algebra of type $A_2$. Thus $\mathfrak{W}(\Lambda_2)$ contains five irreducible morphisms which have source $\add\{P(2),S(2),P(3)\}$ and five morphisms $\add\{P(2),S(2),P(3)\} \rightarrow 0$. The remaining subcategories shown are semisimple, so each is the source of four irreducible morphisms and four morphisms with target 0 in $\mathfrak{W}(\Lambda_2)$.

\begin{center}
\begin{figure}
\begin{tikzpicture}

\begin{scope}[thick,decoration={
	markings,
	mark=at position 0.5 with {\arrow[scale=1.5]{>}}}
	]
        \draw[postaction={decorate}] (0.2,0)--(0.2,4) node [midway,anchor=west]{$P(3)[1]$};
        \draw[postaction={decorate}] (-0.2,0)--(-0.2,4) node [midway,anchor=east]{$P(3)$};
        \draw[postaction={decorate}] (0,0.2)--(3.8,2) node [midway,anchor=south east]{$P(2)$};
        \draw[postaction={decorate}] (0.2,0)--(4,1.8) node [midway,anchor=north west]{$P(2)[1]$};
        \draw[postaction={decorate}] (0,0.2)--(-3.8,2) node [midway,anchor=south west]{$P(1)$};
        \draw[postaction={decorate}] (-0.2,0)--(-4,1.8) node [midway,anchor=north east]{$P(1)[1]$};
        \draw[postaction={decorate}] (0,0)--(5,-2) node [midway,anchor=south west]{$S(1)$};
        \draw[postaction={decorate}] (0,0)--(-5,-2) node [midway,anchor=south east]{$M_{(2,3,0)}$};
        \draw[postaction={decorate}] (0,0)--(3,-4) node [pos=0.7,anchor= west]{$M_{(2,1,0)}$};
        \draw[postaction={decorate}] (0,0)--(-3,-4) node [pos=0.7,anchor=east]{$M_{(3,4,0)}$};
    \end{scope}
    
\draw[dashed] (0,0)--(0,-5);

\node at (0,0) [draw,fill=white] {$\mods\Lambda_2$};
\node at (0,4) [draw,fill=white] {$P(3)^\perp$};
\node at (-4,2) [draw,fill=white] {$\add\{P(2),S(2),P(3)\}$};
\node at (4,2) [draw,fill=white] {$\add\{S(1),P(3)\}$};
\node at (5,-2) [draw,fill=white] {$\add\{M_{(2,1,0)},P(3)\}$};
\node at (3,-4) [draw,fill=white] {$\add\{M_{(3,2,0)},P(3)\}$};
\node at (-5,-2) [draw,fill=white] {$\add\{M_{(1,2,0)},P(3)\}$};
\node at (-3,-4) [draw,fill=white] {$\add\{M_{(2,3,0)},P(3)\}$};

\node at (0.75,-2) {$\cdots$};
\node at (0.75,-3.25) {$\cdots$};
\node at (0.75,-4.5) {$\cdots$};
\node at (-0.75,-2) {$\cdots$};
\node at (-0.75,-3.25) {$\cdots$};
\node at (-0.75,-4.5) {$\cdots$};

\end{tikzpicture}
\caption{The irreducible morphisms in $\mathfrak{W}(\Lambda_2)$ with source $\mods\Lambda_2$. Here, $\Lambda_2 = K(1\rightrightarrows2\rightarrow3)/\mathrm{rad^2}$.}\label{fig:example}
\end{figure}
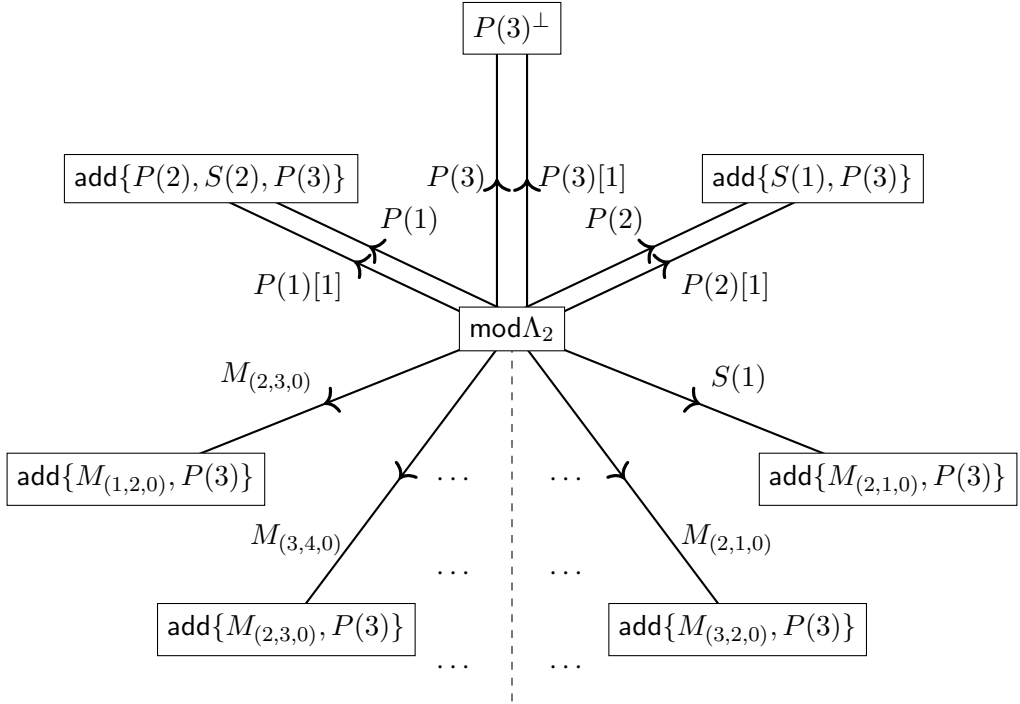
\end{center}

%%%%%%%%%%%%%%%%%%%%%%%%%%%%%%%%%%%%%%

\bibliographystyle{amsalpha}
\bibliography{refs.bib}

\end{document}